\documentclass[10pt,a4paper]{article}

\usepackage{amssymb}
\usepackage{amsmath}
\usepackage{amsthm}
\usepackage{xcolor}
\usepackage{graphicx}
\usepackage{epstopdf}
\usepackage{enumitem}
\usepackage[linesnumbered, boxed]{algorithm2e}
\usepackage{comment}
\usepackage{psfrag}
\usepackage[breaklinks,bookmarks=false]{hyperref}
\usepackage[textsize=tiny,color=yellow!60!green]{todonotes}
\usepackage{multirow}
\usepackage{bm}
\usepackage{mathtools}
\usepackage{stmaryrd}
\usepackage{soul}

\usepackage{mathbbol}
\usepackage{amssymb}             

\DeclareSymbolFontAlphabet{\amsmathbb}{AMSb}%
\newcommand{\defeq}{\vcentcolon=}
\newcommand{\textn}[1]{\textnormal{#1}}

\DeclareMathAlphabet{\mathbfit}{OML}{cmm}{b}{it}
\newcommand{\jump}[1]{\left\llbracket #1 \right\rrbracket}

\newcommand{\RR}{{\amsmathbb R}}

\newcommand{\Omegai}{\Omega_{i}}

\newcommand{\vecv}{{\mathbf{v}}}
\newcommand{\vecr}{{\mathbf{r}}}
\newcommand{\vecu}{{\mathbf{u}}}
\newcommand{\vecw}{{\mathbf{w}}}

\newcommand{\vecF}{{\mathbf{\mathcal{K}}}}

\newcommand{\vK}{\mathbf K}
\newcommand{\II}{{\textbf I}}

\newcommand{\hg}{h,\gamma}

\newcommand{\vecn}{{\mathbf{n}}}

\newcommand{\hi}{h,i}

\newcommand{\eq}{\defeq}

\newcommand{\Tau}{\mathcal{T}}

\DeclareMathOperator*{\argmin}{arg\,min}
\DeclareMathOperator*{\spanvec}{span}
\usepackage{tabularx}

\newtheorem{thm}{Theorem}
\numberwithin{thm}{section}

\newtheorem{rem}[thm]{Remark}

\newtheorem{algo}[thm]{Algorithm}

\newcommand{\bse}{\begin{subequations}}
\newcommand{\ese}{\end{subequations}}

\makeatletter
\newcommand{\eqnum}{\refstepcounter{equation}\textup{\tagform@{\theequation}}}
\makeatother

\graphicspath{{Figures/},{figures/}}

\title{A multiscale flux basis for mortar mixed discretizations of reduced Darcy-Forchheimer fracture models}

\author{Elyes Ahmed\footnotemark[2]
\and Alessio Fumagalli\footnotemark[2]
\and Ana Budi\v{s}a\footnotemark[2]\ 
}
\date{\today}

\begin{document}

\maketitle

\renewcommand{\thefootnote}{\fnsymbol{footnote}}

\footnotetext[2]{Department of Mathematics, University of Bergen, P. O. Box 7800, N-5020 Bergen, Norway.
\href{mailto:elyes.ahmed@uib.no}{elyes.ahmed@uib.no},
\href{mailto:alessio.fumagalli@uib.no}{alessio.fumagalli@uib.no},
\href{mailto:ana.budisa@uib.no}{ana.budisa@uib.no},
}
\renewcommand{\thefootnote}{\arabic{footnote}}

\numberwithin{equation}{section}





\begin{abstract}
    In this paper,  a  multiscale flux basis algorithm  is developed  to
    efficiently solve a flow problem in fractured porous media. Here, we take
    into account a mixed-dimensional setting of the discrete fracture matrix
    model, where the fracture network is represented as lower-dimensional
    object. We assume the linear Darcy model in the rock matrix and the non-linear
    Forchheimer model in the fractures. In our formulation, we are able to
    reformulate the matrix-fracture problem to only the fracture network problem
    and, therefore, significantly reduce the computational cost. The resulting
    problem is then a non-linear interface problem that  can be solved using a
    fixed-point or Newton-Krylov methods, which in each  iteration require
    several solves of Robin  problems in the surrounding rock matrices.       To
    achieve this, the flux exchange (a linear Robin-to-Neumann co-dimensional
    mapping) between the porous medium and the fracture network is done
    \textit{offline} by  pre-computing a multiscale flux basis that consists of
    the flux  response from each  degree of freedom on the fracture network.
    This delivers a \textit{conserve} for the basis that handles the solutions
    in  the rock matrices for each degree of freedom in the fractures pressure
    space.  Then,  any Robin sub-domain problems are  replaced by linear combinations of the multiscale flux basis during the interface iteration. The
    proposed approach is, thus, agnostic to the physical model in the fracture
    network. Numerical experiments demonstrate the   computational gains of
    pre-computing the flux exchange  between the porous medium and the fracture
    network  against  standard non-linear domain decomposition approaches.
\end{abstract}


\vspace{3mm}




\pagestyle{myheadings} \thispagestyle{plain} \markboth{E. Ahmed, A. Fumagalli and A. Budi\v{s}a}{A multiscale flux basis for reduced Darcy-Forchheimer fracture models}




%
\section{Introduction}
\label{sec:Intro}

Using the techniques of domain decomposition~\cite{MR2386967}, a first reduced model has
been proposed for   flow in  a porous medium with  a fracture  in which  the flow
in the fracture is governed by the Darcy-Forchheimer's law while that in the surrounding matrix is
governed by Darcy's law.

We consider here the generalized model given in~\cite{MR3264361}, for which we let $\Omega$ to be a
bounded domain in $\RR^{d}$, $d=2,3$,  with boundary $\Gamma=\partial\Omega$, and we let $\gamma \subset\Omega$ be a $(d-1)$-dimensional surface that
divide  $\Omega$ into two sub-domains: $\Omega=\Omega_{1}\cup\Omega_{2}\cup\gamma$, where $\gamma=\partial \Omega_{1}\cap\partial \Omega_{2}$ and
$\Gamma_{i}=\partial \Omega_{i}\cap \partial \Omega$, $i=1,2$. The reduced  model problem
as presented in~\cite{MR3264361} is  as follows:
\bse\label{Initial_system_porous}\begin{alignat}{5}
\label{Initial_system_d_p}\vK^{-1}_{i}\vecu_{i}+\nabla p_{i}  &= 0 &&\quad \textn{in} \;  \Omega_{i},\\
\label{Initial_system_c_p}\nabla\cdot\vecu_{i} &=f_{i} &&\quad \textn{in} \;
\Omega_{i},\\
\label{Initial_system_bd_p}p_{i} &=0 &&\quad \textn{in} \;   \Gamma_{i},
\end{alignat}\ese
for $i=1,2$, together with
\bse\label{Initial_system_fracture}\begin{alignat}{4}
\label{Initial_system_d_f}(\vK^{-1}_{\gamma}+\beta_{\gamma}\II|\vecu_{\gamma}|)\vecu_{\gamma}  &= -\nabla_{\tau} p_{\gamma} &&\quad \textn{in} \;  \gamma,&\\
\label{Initial_system_c_f}\nabla_{\tau}\cdot\vecu_{\gamma} &=f_{\gamma}+\left(\vecu_{1}\cdot\vecn_{1}+\vecu_{2}\cdot\vecn_{2}\right) &&\quad \textn{in} \;  \gamma,&\\
\label{Initial_system_bd_f} p_{\gamma} &=0 &&\quad \textn{in} \;   \partial\gamma,&
\end{alignat}\ese
and subject to  the following  interface  conditions
\begin{alignat}{4}
\label{Initial_system_r_interface}-\vecu_{i}\cdot\vecn_{i}+\alpha_{\gamma} p_{i} &=\alpha_{\gamma} p_{\gamma} &&\quad \textn{on} \;   \gamma,
\end{alignat}
for $i=1,2$. Here, $\nabla_{\tau}$ denotes  the $(d-1)$-dimensional   gradient  operator  in  the  plane  of $\gamma$, the coefficient
$\vK_{i}$, $i=1,2$, is the hydraulic conductivity tensor in the sub-domain
$\Omegai$, and $\vK_{\gamma}$  is the hydraulic conductivity tensor in the fracture,
$\II\in\RR^{d\times d}$ is the identity matrix, $\vecn_{i}$ is the outward unit normal vector to $\partial\Omegai$,
and  $\beta_{\gamma}$ is a  non-negative scalar known  as the Forchheimer
coefficient.
\begin{figure}[htb]
    \centering
    \resizebox{0.2\textwidth}{!}{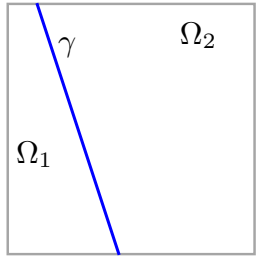}%
    \caption{Graphical example of problem
    \eqref{Initial_system_porous}-\eqref{Initial_system_r_interface}.}%
    \label{fig:single_domain}
\end{figure}
In~\eqref{Initial_system_r_interface}, the coefficient $\alpha_{\gamma}$ is a
function  proportional to  the normal component of the permeability of the
physical fracture and inversely proportional to the fracture width/aperture. We
refer to~\cite{Frih2008} for a more detailed model description.  For
illustration purposes, we give a simple graphical example of a fractured porous
medium in Figure \ref{fig:single_domain}.

The system \eqref{Initial_system_porous}-\eqref{Initial_system_r_interface} can be seen as a domain decomposition problem,
with non-standard and non-local boundary conditions between the sub-domains $\Omegai$,
$i=1,2$. The equations~\eqref{Initial_system_porous} are the mass conservation equation and the Darcy's law equation in the
sub-domain $\Omegai$ while equations~\eqref{Initial_system_fracture} are the lower-dimensional mass conservation and the Darcy-Forchheimer equation in the fracture of co-dimension
$1$. The  last equation~\eqref{Initial_system_r_interface} can be seen as a
Robin boundary condition   for the sub-domain $\Omega_{i}$
with a dependence on the pressure on the fracture $\gamma$.
Clearly, if $\beta_{\gamma}=0$, then~\eqref{Initial_system_fracture} is reduced  to  a linear Darcy flow in the fracture.
The homogeneous Dirichlet boundary conditions~\eqref{Initial_system_bd_p}
and~\eqref{Initial_system_bd_f}
are considered  merely for simplicity. The functions $f_{i}\in L^{2}(\Omegai)$, $i=1,2$
and $f_{\gamma}\in L^{2}(\gamma)$ are source terms in the matrix and in the fracture, respectively.

The mixed-dimensional problem~\eqref{Initial_system_porous}-\eqref{Initial_system_r_interface} is an  alternative
to the possibility of using a very fine grid in the physical fracture and a necessarily much coarser grid
away from the fracture. This idea was developed in~\cite{alboin1999} for highly permeable fractures and
in~\cite{MR2512496}  for fractures that may be highly permeable or nearly impermeable.
We  also  refer  to~\cite{MR2776916,MR2585104,MR3489127}  for  similar  models.  For  all  of  the  above
models,  where the  linear  Darcy's  law  is  used  as  the  constitutive  law  for  flow  in  the fractures  as  well  as  in  the  surrounding
domains,   there are interactions between fractures and surrounding domains.  This coupling is ensured
using    Robin type  conditions  as in~\cite{MR2142590},    delivering   discontinuous normal velocity and pressure
across the fractures.  Particularly,   for  fractures  with  large  enough  permeability,  Darcy's  law
is   replaced  by  Darcy-Forchheimer's law as established in~\cite{MR3264361}, which  complicates the coupling with the surrounding medium.

Several numerical schemes have been developed for fracture models, such as a cell-centered finite volume scheme in~\cite{HAEGLAND20091740},
an  extended finite element method  in~\cite{MR3631391}, a mimetic finite difference~\cite{MR3507274} and a  block-centered finite difference method
in~\cite{MR3693346}. The  aforementioned  numerical  approaches     solve  coupled  fracture models  directly.
However,  different equations defined in different regions are
varied in type,  such as coupling linear and non-linear systems, and often interface conditions involve new   variables in  different  domains,
which results  in  very  complex  algebraic structures. Particularly,  several papers deal with the analysis and implementation of mixed methods
applied to the above model problem in the linear  case,   on conforming and
non-conforming
grids~\cite{MR2142590,angelo_scotti_2012,Frih2012,MR2773340,MR2890281}. In~\cite{MR2386967},
the  model problem~\eqref{Initial_system_porous}-\eqref{Initial_system_r_interface}  was solved using domain decomposition techniques based on
mixed finite element methods (see~\cite{alboin1999} for the linear counterpart).

The purpose of this paper is to  propose an efficient domain decomposition
method to
solve~\eqref{Initial_system_porous}-\eqref{Initial_system_r_interface} based on
the multiscale mortar mixed finite element method (MMMFEM) \cite{Ganis2009}. The method
reformulates~\eqref{Initial_system_porous}-\eqref{Initial_system_r_interface}
into   an  interface problem by eliminating the  sub-domain variables. The
resulting interface problem is a superposition of a non-linear operator handling
the flow  on the fracture and a linear operator presenting the  flux
contribution from the sub-domains.  When applying the MMMFEM,  an \textit{outer--inner iterative algorithm}  like,  the  Newton--GMRes
(or any Krylov solver) method or  fixed-point--GMRes method,  is used to solve the interface
problem. As an example,  if a fixed-point  method (outer) is adopted, the linearized
interface equation for the interface update can  be solved with a domain
decomposition algorithm (inner),    in which  at each
iteration sub-domain solves,  together with  inter-processor communication, are
required.  The main issue of this outer-inner algorithm is that it leads to an
excessive calculation  from the sub-domains, as the dominant computational cost is measured by the number of sub-domain solves.

The new implementation recasts this algorithm by distinguishing  the linear and non-linear contributions in the
overall calculation and employing the multiscale flux basis functions
from~\cite{MR2557486} for the linear part of the problem,  before the non-linear
interface iterations begin. The fact that the non-linearity
in~\eqref{Initial_system_porous}-\eqref{Initial_system_r_interface}  is only
within the fracture, we can adopt the notion that sub-domain problems can be
expressed  as a superposition of multiscale basis functions. In our terminology
the mortar variable considered in \cite{Ganis2009} becomes the fracture pressure,
these multiscale
flux basis with respect to the fracture pressure can be
computed by  solving a fixed number of Robin sub-domain problems, that is equal
to the number of fracture pressure degrees of freedom per sub-domain. Furthermore,
this is done in parallel without any inter-processor communication.

An inexpensive linear combination of the multiscale flux basis functions  then
circumvents the need to solve any sub-domain problems in the inner domain
decomposition iterations. This procedure can be enhanced by applying  interface
preconditioners as in~\cite{amir2006decomposition,MR2142590,MR3457700} and by
using a posteriori error estimates
of~\cite{Pench_Voh_Whee_Wild_a_post_MS_MN_M_13} to adaptively refine the  mesh
grids. This calculation made in an \textit{offline} step   typically spares
numerous  unnecessary sub-domain solves. Precisely, in the original
implementations, the number of sub-domain solves   is approximately equal to
$\sum_{k=1}^{N_{\textnormal{lin}}} N^{i}_{\textnormal{dd}}$, where
$N_{\textnormal{lin}}$ is the number of  iterations of the linearization
procedure, and  $N^{k}_{\textnormal{dd}}$ denotes the number of domain
decomposition iterations (GMRes or any Krylov solver). For the new
implementation,  the number of sub-domains solves will be reduced if
$\sum_{k=1}^{N_{\textnormal{lin}}} N^{k}_{\textnormal{dd}}$ exceeds the maximum
number of fracture pressure degrees of freedom on any sub-domain.

This step of \textit{freezing}  the contributions on the flow from the rock
matrices can be  \textit{easily coded,  cheaply evaluated},  and
\textit{efficiently used} in  \textit{practical simulations}, i.e,  it permits
reusing the same basis functions  to
extend~\eqref{Initial_system_porous}-\eqref{Initial_system_r_interface} to
simulate various linear and non-linear models for flow in the fracture, such
as generalized Forchheimer's laws:
\begin{gather*}
    (\vK^{-1}_{\gamma}+\beta_{\gamma}\II|\vecu_{\gamma}|+\zeta_{\gamma}\II|\vecu_{\gamma}|^{2})\vecu_{\gamma}
    = -\nabla_{\tau} p_{\gamma},\\
    (\vK^{-1}_{\gamma}e^{\zeta
    p_\gamma}+\beta_{\gamma}\II|\vecu_{\gamma}|)\vecu_{\gamma}  = -\nabla_{\tau} p_{\gamma},
\end{gather*}
as well as exploring  the  fracture and barrier cases and comparing in a cheap
way various non-linear solvers to~\eqref{Initial_system_porous}-\eqref{Initial_system_r_interface}. Crucially, the present approach
can naturally be integrated into  discrete fracture networks (DFNs)
models~\cite{MR2773340,MR2890281,MR3757111,MR3264347}, which in contrast to
discrete fracture models (DFMs), do not consider the flow in the surrounding sub-domains, but
handle  both a large number of fractures and a complex interconnecting network
of these fractures.


For the presenting setting,
we allow for the discretization of~\eqref{Initial_system_porous}-\eqref{Initial_system_r_interface} by different numerical methods applied separately
in the surrounding sub-domains and in the fracture. We allow for the cases where the grids of the porous
sub-domains do not match along the fracture, where different mortar grid elements are used. We also investigate the case
where the permeability in the fracture $\vK_{\gamma}$ is much lower than the permeability in the surrounding matrix $\vK$.

The library PorePy \cite{Keilegavlen2017a} has been used and extended to cover
the numerical schemes and examples introduced in this article. The main
contribution to the library is the implementation of the multiscale and domain
decomposition frameworks. Even if we focus on lowest-order Raviart-Thomas-N\'{e}d\'{e}lec
finite elements, our implementation is agnostic with respect to the numerical
scheme. The example presented are also available in the GitHub
repository.

This paper is organized as follows:  Firstly, the variational formulation of the
problem and the MMMFEM approximation are given in
Section~\ref{sec:problem_formulation}. Therefrom, the reduction of the original
problem into nonlinear interface problem is introduced.  The
linearization--domain-decomposition procedures are formulated in
Section~\ref{sec:nonlinear_algorithms}. Section~\ref{sec:Mufbi} describes the
implementation based on the multiscale flux basis. We show  that structurally
the same implementation can  be extended for more complex intersecting fractures
model. Finally, we showcase the performance of our method on several numerical
examples in Section~\ref{sec:examples} and draw the conclusions in
Section~\ref{sec:conclusion}.
%

\section{Non-linear domain decomposition method}
\label{sec:problem_formulation}
 As explained earlier, it is natural to solve the mixed-dimensional problem~\eqref{Initial_system_porous}-\eqref{Initial_system_r_interface} using domain decomposition techniques.
 To this aim,   we introduce the weak spaces in each sub-domain $\Omega_{i}$, $i=1,2,$
 \begin{alignat*}{2}
&\mathbf{V}_{i}\defeq\mathbf{H}(\textn{div}, \Omega_{i}),\qquad &&  M_{i}\defeq L^{2}(\Omega_{i}),
\end{alignat*}
and   define their global versions by
\begin{alignat*}{2}
&\mathbf{V}\defeq\bigoplus_{i=1}^{2}\mathbf{V}_{i},\qquad&&  M\defeq\bigoplus_{i=1}^{2} M_{i}.
\end{alignat*}
Equivalently, we introduce the weak spaces on the fracture $\gamma$, i.e,
\begin{alignat*}{2}
\mathbf{V}_{\gamma}\defeq\mathbf{H}(\textn{div}_{\tau}, \gamma),\qquad M_{\gamma}\defeq L^{2}(\gamma).
\end{alignat*}
Following~\cite{MR2386967,MR2142590}, a mixed-dimensional weak form of~\eqref{Initial_system_porous}-\eqref{Initial_system_r_interface} asks for $(\vecu,p)\in \mathbf{V}\times M$ and
$(\vecu_{\gamma},p_{\gamma})\in \mathbf{V}_{\gamma}\times M_{\gamma}$
such that, for each $i \in \{1,2\}$,
\bse \label{weak_mixed_formulation} \begin{alignat}{2}
\nonumber&(\vK^{-1}\vecu,\vecv)_{\Omega_{i}}+\alpha_{\gamma}^{-1}\langle
\vecu\cdot\vecn_{i},\vecv\cdot\vecn_{i}\rangle_{\gamma}
&&\\
&\qquad\quad\qquad\quad=(p,\nabla\cdot\vecv)_{\Omega_{i}}-\langle
p_{\gamma}, \vecv\cdot\vecn_{i} \rangle_\gamma\quad&&\forall\vecv\in \mathbf{V},\\
&(\nabla\cdot \vecu,q)_{\Omega_{i}}=(f,q)_{\Omega_{i}} \quad&&\forall q\in  M,\\
&\langle \vecF^{-1}(\vecu_{\gamma})\vecu_{\gamma}
,\vecv_{\gamma}\rangle_\gamma=\langle p_{\gamma},
\nabla_{\tau}\cdot\vecv_{\gamma} \rangle_\gamma\quad&& \forall\vecv_{\gamma}\in \mathbf{V}_{\gamma},\\
&\langle \nabla_{\tau} \cdot\vecu_{\gamma}, q_{\gamma} \rangle_\gamma =
\langle f_{\gamma} + \jump{\vecu \cdot \vecn}, q_{\gamma} \rangle_\gamma \quad && \forall q_{\gamma}\in M_{\gamma},
\end{alignat}
\ese
where we introduced  the functions $\vK$ and $f$  in
 $\Omega_{1}\cup \Omega_{2}$  such that $\vK_{i}=\vK|_{\Omegai}$, and
 $f_{i}=f|_{\Omegai}$, $i=1,2$. The jump $\jump{\cdot}$
 is  defined by
 \begin{equation*}
  \jump{\vecu \cdot \vecn}\defeq\vecu_{1}\cdot\vecn_{1}+\vecu_{2}\cdot\vecn_{2},
 \end{equation*}
with $\vecn_{i}$ the outer unit normal vector of $\Omega_i$ on $\gamma$, for $ i=1,2 $. Finally, the non-linear
term is defined as
\begin{gather*}
    \vecF^{-1}(\vecu_{\gamma}) \defeq \vK^{-1}_{\gamma}+\beta_{\gamma}\II|\vecu_{\gamma}|,
\end{gather*}
The reader is referred to~\cite{MR3264361} 
for proof of the existence and uniqueness of a solution to the variational formulation~\eqref{weak_mixed_formulation}. 


\subsection{The discrete problem}
\label{sec:problem_form.subsec:mfe}

Let $\Tau_{\hi}$ be a partition of the sub-domain
$\Omegai$ into either $d$-dimensional simplicial and/or
rectangular elements. We also
let  $\Tau_{\hg}$ to be a partition of the fracture $\gamma$ into
$(d-1)$-dimensional simplicial and/or  rectangular elements. Note that, for both partitions,
general  elements can be treated via sub-meshes, see~\cite{MR3033022} and the references therein. Moreover, we assume that each  partition is conforming
within each sub-domain as well as in the fracture. The meshes $\Tau_{\hi}$,
$i=1,2$, are  allowed to be non-conforming on the
fracture-interface $\gamma$, but also different
from  $\Tau_{\hg}$. We then set $\Tau_{h}=\displaystyle{\cup_{i=1}^{2}\Tau_{\hi}}$ and denote by $h$ the
maximal element diameter in~$\Tau_{h}$.
\begin{figure}[htb]
    \centering
    \includegraphics[width=0.45\textwidth]{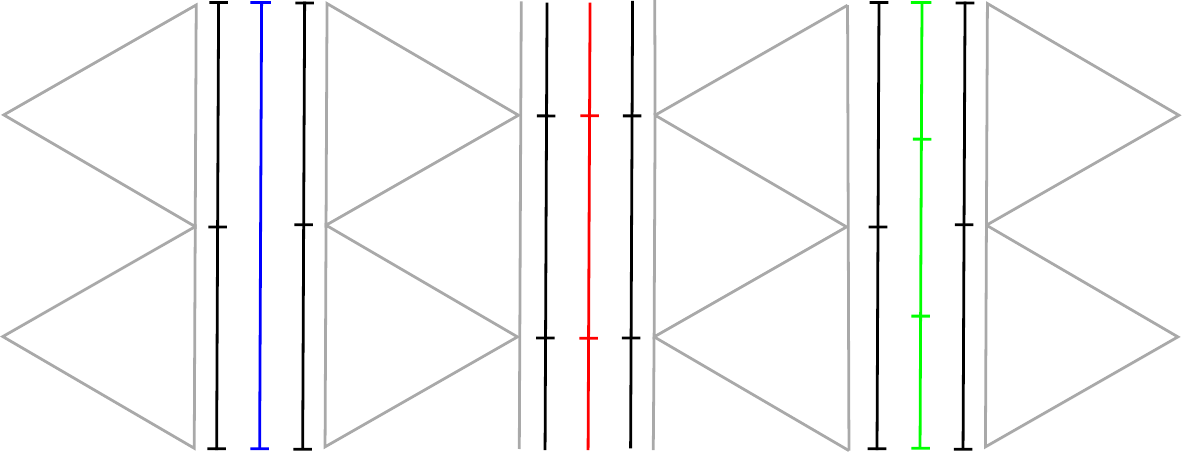}%
    \caption{Representation of three possible fracture mesh configurations: on
    the left coarser, on the centre conforming, and on the right finer. The
    triangles are represented in grey.}
    \label{fig:mortar}
\end{figure}
For  the scalar unknowns, we
introduce the   approximation spaces
$M_{h}=M_{h,1}\times M_{h,2}$ and $M_{\hg}$, where $M_{\hi}$, $i=1,2$, respectively
$M_{\hg}$, is the space of piecewise constant functions associated with $\Tau_{\hi}$,
$i=1,2$, respectively $\Tau_{\hg}$. For  the vector unknowns, we
introduce the   approximation spaces
$\mathbf{V}_{h}=\mathbf{V}_{h,1}\times \mathbf{V}_{h,2}$ and $\mathbf{V}_{\hg}$, where
$\mathbf{V}_{\hi}$, $i=1,2$ and
$\mathbf{V}_{\hg}$, are the lowest-order Raviart-Thomas-N\'{e}d\'{e}lec finite elements spaces
associated with  $\Tau_{\hi}$,  $i=1,2$ and $\Tau_{\hg}$, respectively. Clearly, in contrast to what is done in \cite{Ganis2009,Efendiev2009},
we use the same order of the polynomials for the interface-pressure and the normal traces of the sub-domain velocities on the interface.

The  discrete  mixed-dimensional  finite  element  approximation
of~\eqref{weak_mixed_formulation}  is  as
follows: find $(\vecu_{h},p_{h})\in
\mathbf{V}_{h}\times M_{h}$ and
$(\vecu_{h,\gamma},p_{h,\gamma})\in \mathbf{V}_{h,\gamma}\times M_{h,\gamma}$
such that, for each $i\in\{1,2\}$,
\bse\label{discrete_mixed_formulation}\begin{alignat}{2}
    \nonumber&(\vK^{-1}\vecu_{h},\vecv)_{\Omega_{i}}+\alpha_{\gamma}^{-1}\langle
    \vecu_{h}\cdot\vecn_{i},\vecv\cdot\vecn_{i} \rangle_{\gamma}
    &&\\
    &\qquad\quad\qquad\quad=(p_{h},\nabla\cdot\vecv)_{\Omega_{i}}-\langle p_{\hg},
    \vecv\cdot\vecn_{i} \rangle_\gamma\quad&&\forall\vecv\in \mathbf{V}_{h},\\
    &(\nabla\cdot \vecu_{h},q)_{\Omega_{i}}=(f,q)_{\Omega_{i}} \quad&&\forall q\in  M_{h},\\
&\langle \vecF^{-1}(\vecu_{\hg})\vecu_{\hg},
\vecv_{\gamma}\rangle_\gamma=\langle p_{\hg}, \nabla_{\tau}\cdot\vecv_{\gamma}
\rangle_\gamma\quad&& \forall\vecv_{\gamma}\in \mathbf{V}_{\hg},\\
&\langle \nabla_{\tau} \cdot\vecu_{\hg}, q_{\gamma} \rangle_\gamma=
\langle f_{\gamma}+\jump{\vecu_{h}\cdot\vecn}, q_{\gamma}\rangle_\gamma\quad&& \forall q_{\gamma}\in M_{\hg}.
\end{alignat}
\ese

The next step in formulating   a  multiscale flux basis algorithm   to solve~\eqref{discrete_mixed_formulation} is to adopt  domain decomposition techniques
to reduce the global mixed-dimensional problem to
an interface problem posed only on the fracture~\cite{Amir2005}.

\subsection{Reduction to interface problem}
\label{sec:reduced_problem}
We introduce the discrete (linear) Robin-to-Neumann  operator $\mathcal{S}^{\textn{RtN}}_{i}$, $i=1,2$:
\begin{alignat*}{2}
\mathcal{S}^{\textn{RtN}}_{i}\,\colon\,L^{2}(\gamma)\times L^{2}(\Omegai)&\rightarrow L^{2}(\gamma),\\
 \mathcal{S}^{\textn{RtN}}_{i}(\lambda,f)&=-\vecu_{h}(\lambda,f)\cdot\vecn_{i},
\end{alignat*}
where
$(\vecu_{h},p_{h})\in
\mathbf{V}_{h}\times M_{h}$ is the solution of the sub-domain problems with source term $f$, homogeneous
Dirichlet boundary condition on $ \partial\Omega$, and $\lambda$ as a Robin boundary condition along the fracture $\gamma$, i.e, for $i=1,2$,
\bse\label{subdomain_problems}\begin{alignat}{2}
\nonumber&(\vK^{-1}\vecu_{h},\vecv)_{\Omega_{i}}+\alpha_{\gamma}^{-1}\langle
\vecu_{h}\cdot\vecn_{i},\vecv\cdot\vecn_{i} \rangle_{\Omega_{i}}\\
&\qquad\quad\qquad\quad=(p_{h},\nabla\cdot\vecv)_{\Omega_{i}} -\langle \lambda,
\vecv\cdot\vecn_{i}\rangle_\gamma\quad&& \forall \vecv\in\mathbf{V}_{h,i},\\
&(\nabla\cdot \vecu_{h},q)_{\Omega_{i}}=(f,\vecv)_{\Omega_{i}}&&\forall q\in M_{h,i}.
\end{alignat}
\ese
Then we set
\begin{alignat*}{2}
\mathcal{S}^{\textn{RtN}}\,\colon\,L^{2}(\gamma)\times L^{2}(\Omega)&\rightarrow L^{2}(\gamma),\\
 \mathcal{S}^{\textn{RtN}}(\lambda,f)&=\sum_{i=1}^{2}\mathcal{S}^{\textn{RtN}}_{i}(\lambda,f_{i}).
\end{alignat*}
With these notations, we can see that solving~\eqref{discrete_mixed_formulation} is equivalent to solving the following non-linear mixed
interface problem: find
$(\vecu_{h,\gamma},p_{h,\gamma})\in \mathbf{V}_{h,\gamma}\times M_{h,\gamma}$
such that,
\bse\label{discrete_mixed_formulation_interface}\begin{alignat}{2}
&\langle \vecF^{-1}(\vecu_{h,\gamma})\vecu_{h,\gamma},
\vecv_{\gamma}\rangle_\gamma-\langle p_{\hg}, \nabla_{\tau}\cdot\vecv_{\gamma}
\rangle_\gamma=0\quad&& \forall\vecv_{\gamma} \in \mathbf{V}_{h,\gamma},\\
&\langle \nabla_{\tau} \cdot\vecu_{h,\gamma}, q_{\gamma}\rangle_\gamma+\langle
\mathcal{S}^{\textn{RtN}}(p_{\hg},f), q_{\gamma}\rangle_\gamma=\langle
f_{\gamma}, q_{\gamma}\rangle_\gamma\quad&& \forall q_{\gamma}\in M_{h,\gamma},
\end{alignat}
\ese
or equivalently
\bse\label{discrete_mixed_formulation_interface2}\begin{alignat}{2}
&\langle \vecF^{-1}(\vecu_{h,\gamma})\vecu_{h,\gamma},
\vecv_{\gamma}\rangle_\gamma-\langle p_{\hg}, \nabla_{\tau}\cdot\vecv_{\gamma}
\rangle_\gamma=0\quad&& \forall\vecv_{\gamma} \in \mathbf{V}_{h,\gamma},\\
&\langle \nabla_{\tau} \cdot\vecu_{h,\gamma}, q_{\gamma}\rangle_\gamma+\langle
\mathcal{S}_{\gamma}(p_{\hg}), q_{\gamma}\rangle_\gamma=\langle
f_{\gamma}+g_{\gamma}, q_{\gamma}\rangle_\gamma\quad&& \forall q_{\gamma}\in M_{h,\gamma},
\end{alignat}
\ese
where we  have set
\begin{alignat}{2}\label{right_hand_side}
 &\mathcal{S}_{\gamma}(\lambda)=\mathcal{S}^{\textn{RtN}}(\lambda,0)\quad\text{and}\quad g_{\gamma}=-\mathcal{S}^{\textn{RtN}}(0,f).
\end{alignat}
The  above distinction  is classical in domain decomposition techniques
in which we split the sub-domain problems  into two families of local problems on each $\Omegai$: one is with zero source and specified  Robin value
on the fracture-interface, and the other is with zero  Robin value on the fracture-interface and specified source.
In compact form, the mixed interface Darcy-Forchheimer problem~\eqref{discrete_mixed_formulation_interface2} can be rewritten as
\begin{align}\label{compact_system}
\begin{bmatrix}
 \vecF^{-1}(\cdot)& B_{\gamma}^{\textn{T}}\\[2mm]
 B_{\gamma}& \mathcal{S}_{\gamma}\\[2mm]
\end{bmatrix}\begin{bmatrix}
 \vecu_{\gamma}\\[2mm]
 p_{\gamma}\\[2mm]
 \end{bmatrix} =  \begin{bmatrix}
 0\\[2mm]
 g_{\gamma}+f_{\gamma}
\end{bmatrix}.
\end{align}

This  system is a non-linear mixed interface problem~\cite{MR3392446} that can  be solved iteratively by using fixed
point iterations or via a Newton-Krylov method.  To present the two approaches, let us first consider  the  linear context, i.e,
suppose the operator $\vecF^{-1}(\cdot)$ is linear.  Then~\eqref{compact_system} is the  system associated to the
linear mixed Darcy problem on the fracture that can be solved  using a Krylov type method, such as GMRes or MINRes. Given an initial
guess $\vecw_{\gamma}^{(0)}=[\vecu_{\gamma}^{(0)},p^{(0)}_{\gamma}]^{\textn{T}}$, the GMRes
algorithm computes
\begin{alignat}{2}\label{gmres_pseudocode}
 \vecw_{\gamma}^{(m)}=\underset{\vecv \, \in \,
 \vecw_{\gamma}^{(0)}+\mathcal{K}_{m}(\mathcal{A}_{\gamma},
 \vecr^{(0)}_{\gamma})}{\argmin}||\mathbf{b}_{\gamma}-\mathcal{A}_{\gamma}\vecv||_{2}
 \quad \text{for} \quad m\geq 1,
\end{alignat}
as an approximate solution to~\eqref{compact_system}, where  $\mathcal{A}_{\gamma}$ is the associated stiffness matrix of the linear system,
$\mathbf{b}_{\gamma}$ is the right-hand side, and
$\mathcal{K}_{m}(\mathcal{A}_{\gamma}, \vecr_{\gamma}^{(0)})$ is the $m$-dimensional
Krylov subspace generated by the initial residual $\vecr_{\gamma}^{(0)}=\mathbf{b}_{\gamma}-\mathcal{A}_{\gamma}\vecw_{\gamma}^{(0)}$, i.e,
\begin{alignat*}{2}
 \mathcal{K}_{m}(\mathcal{A}_{\gamma}, \vecr_{\gamma}^{(0)})=\spanvec (
 \vecr_{\gamma}^{(0)}, \mathcal{A}_{\gamma}\vecr^{(0)},\cdots,
 \mathcal{A}_{\gamma}^{(m-1)}\vecr_{\gamma}^{(0)} ).
\end{alignat*}
Clearly,
each GMRes iteration needs to evaluate
the action of the Robin-to-Neumann type operator $ \mathcal{S}_{\gamma}$ via \eqref{right_hand_side}, representing physically the contributions on the flow from the rock matrices, i.e, to
solve one Robin sub-domain problem per sub-domain. Thus the GMRes algorithm is
implemented in the matrix-free context~\cite{MR3022024,MR2557486,MR2032407}.

One can easily observe that  the evaluation
of  $\mathcal{S}_{\gamma}$  dominates the total computational costs in~\eqref{gmres_pseudocode}.
In practice, this step  is done in parallel  and involves inter-processor communication across the fracture-interface~\cite{MR3022024}.  To present the
evaluating algorithm of  $ \mathcal{S}_{\gamma}$,
we  let $\mathcal{D}^{\textn{T}}_{h,i}\colon \mathbf{V}_{h,i}\cdot\vecn_{i}|_{\gamma}\rightarrow M_{\hg}$ be the $L^{2}$-orthogonal projection  from the
normal trace of the velocity space onto the mortar space
normal trace of the velocity space in sub-domain $\Omega_{i}$, $i=1,2$,
onto the pressure space on the fracture $M_{\hg}$. We then summarize the evaluation of
the interface operator  by the following steps:
\begin{algo}[Evaluating  the action of $ \mathcal{S}_{\gamma}$]\label{Eval_op}~
{
\setlist[enumerate]{topsep=0pt,itemsep=-1ex,partopsep=1ex,parsep=1ex,leftmargin=1.5\parindent,font=\upshape}
\begin{enumerate}
    \item Enter an interface data $\varphi$.
    \item \textn{\textbf{For}} $i=1:2$
    \begin{enumerate}
        \item Project mortar data onto sub-domain boundary, i.e,
        $$ \varphi \overset{\mathcal{D}_{h,i}}{\longrightarrow}\lambda.$$
        \item Solve the sub-domain problem~\eqref{subdomain_problems}
        with Robin boundary condition~$\lambda$ and with $f=0$.
        \item Project  the resulting  flux onto  the mortar space $M_{\hg}$, i.e,
        $$ -\vecu_{h}(\lambda,0)\cdot\vecn_{i}\overset{\mathcal{D}^{\textn{T}}_{h,i}}{\longrightarrow}
        -\mathcal{D}^{\textn{T}}_{h,i}\vecu_{h}(\lambda,0)\cdot\vecn_{i}.$$
    \end{enumerate}
    \textn{\textbf{EndFor}}
    \item Compute the flow contribution from the sub-domains to the
    fracture given by the flux jump across the fracture, i.e,
    $$\mathcal{S}_{\gamma}(\varphi)=\displaystyle\sum_{i=1}^{2}-\mathcal{D}^{\textn{T}}_{h,i}\vecu_{h}(\lambda,0)\cdot\vecn_{i}.$$
\end{enumerate}
}
\end{algo}

\section{Non-linear interface iterations} \label{sec:nonlinear_algorithms}
In this section,  we form two linearization--domain-decomposition algorithms  to solve  the mixed interface Darcy-Forchheimer problem~\eqref{discrete_mixed_formulation_interface2}. For the linearization (outer) of~\eqref{discrete_mixed_formulation_interface2}, a
 first algorithm is based on a  fixed-point method is presented  and the second one is based on  Newton-GMRes method~\cite{MR2597810,MR2948707}. For the solver of the inner systems  (domain decomposition systems), both  methods uses the GMRes method~\eqref{gmres_pseudocode} to solve the reduced mixed interface problems. Note that the two approaches   have competitive  performance
for such nonlinear model problems and  they lead to
different applications of the multiscale flux basis functions of Section~\ref{sec:Mufbi}.
\subsection{Method~1: fixed-point-GMRes}\label{subsec:fixed_point}
  We consider first a standard fixed-point approach  to solve the interface Darcy-Forchheimer problem~\eqref{discrete_mixed_formulation_interface2} (see~\cite{MR2948707}).
Given an initial value $\vecu^{(0)}_{\hg}$, being the solution  of a linear Darcy, for $k=1,2,\dots,\textnormal{until convergence}$, find $(\vecu^{(k)}_{h,\gamma},p^{(k)}_{h,\gamma})\in \mathbf{V}_{h,\gamma}\times M_{h,\gamma}$ such that,
\bse\label{discrete_mixed_formulation_iterative}\begin{alignat}{2}
&\langle \vecF^{-1}(\vecu^{(k-1)}_{h,\gamma}) \vecu^{(k)}_{h,\gamma}, \vecv_{\gamma}
\rangle_\gamma - \langle p^{(k)}_{\hg}, \nabla_{\tau}\cdot\vecv_{\gamma}
\rangle_\gamma=0\quad&&\forall \vecv_{\gamma}\in \mathbf{V}_{h,\gamma},\\
&\langle \nabla_{\tau} \cdot \vecu^{(k)}_{h,\gamma}, q_{\gamma} \rangle_\gamma +
\langle \mathcal{S}_{\gamma}(p^{(k)}_{\hg}), q_{\gamma} \rangle_\gamma =
\langle f_{\gamma}+g_{\gamma}, q_{\gamma}\rangle_\gamma\quad &&\forall q_{\gamma}\in M_{h,\gamma}.
\end{alignat}
\ese
This process is  linear and can be solved using GMRes method~\eqref{gmres_pseudocode}, where  each  iteration needs to set up  the action of the Robin-to-Neumann operator
 $\mathcal{S}_{\gamma}$ using Algorithm~\ref{Eval_op}.
%
The above  fixed-point-GMRes algorithm is iterated  until a fixed-point residual tolerance reaches some prescribed value.

The  result of this procedure is then used to generate
the solution in the sub-domains via
\bse\label{subdomain_solutions}\begin{align}
\vecu_{h}|_{\Omegai}=\vecu_{h}(p^{(\infty)}_{\hg},0)|_{\Omegai}+\vecu_{h}(0,f_{i}),\\
p_{h}|_{\Omegai}=p_{h}(p^{(\infty)}_{\hg},0)|_{\Omegai}+p_{h}(0,f_{i}),
\end{align}\ese
for $i=1,2$, requiring  two additional  sub-domain solves, and where
$p^{(\infty)}_{\hg}$ indicates the fracture pressure at convergence.
\begin{rem}[An alternative to~\eqref{discrete_mixed_formulation_iterative}]
A well-known drawback of  GMRes algorithm for solving  the interface-fracture
problem~\eqref{discrete_mixed_formulation_iterative} is that the number of
iterations depends essentially on the number of sub-domain solves. A preconditioner is then necessary to reduce the
iterations number to a reasonable level. To this aim, it is
possible to
reformulate~\eqref{discrete_mixed_formulation_iterative}  into a
primal problem:   at the iteration $k\geq 1$, by
solving
 for the sole  scalar unknown  $p^{(k)}_{\hg}$, such that
 \bse\label{discrete_primal_formulation_iterative}\begin{align}
  -\nabla_{\tau} \cdot [ -\vecF(p^{(k-1)})\nabla_{\tau} p^{(k)}_{\hg}
  ]+\mathcal{S}_{\gamma}(p^{(k)}_{\hg})=g_{\gamma} +f_{\gamma} \quad\text{on
  }\gamma,
 \end{align}\ese
which can be discretized with a  cell-centered finite volume method, leading to
a symmetric  and positive definite system that can be solved with a CG method.
The CG method can be equipped with a preconditioner being the  inverse of the
discrete counterpart of the operator
$-\nabla_{\tau} \cdot [ -\vecF(p^{(k-1)})\nabla_{\tau}
  ]$ (see~\cite{amir2006decomposition,MR3624734} for more details).
\end{rem}
\begin{rem}[The total computational costs]\label{remark:FVCGscheme}
The total computational costs in the  inner-outer iterative approach~\eqref{discrete_mixed_formulation_iterative} is dominated by
the number  of sub-domain solves required. Precisely, the
total number of sub-domain solves
is given by  $\sum_{k=1}^{N_{\textnormal{lin}}} N^{k}_{\textnormal{dd}}$, where $N_{\textnormal{lin}}$ is the number of  iterations of the fixed-point procedure as outer-loop algorithm, and  $N^{k}_{\textnormal{dd}}$ denotes the number of inner loop domain decomposition iterations (GMRes) at the fixed-point iteration $k\geq1$.
\end{rem}

\subsection{Method~2: Newton-GMRes}\label{subsec:Newton-Gmres}
\label{sec:interface_problem.subsec:newton}
 In the second approach,  we propose Newton's method to solve
the interface Darcy-Forchheimer problem~\eqref{discrete_mixed_formulation_interface2}. For simplicity of notation, we introduce the following
\bse\begin{align*}
  &\vecF^{-1, (k)}\defeq\vecF^{-1}(\vecu^{(k)}_{\hg})=\vK^{-1}_{\gamma}+\beta_{\gamma}\II|\vecu^{(k)}_{h,\gamma}|,\\
  &\vecF_{\partial}^{-1, (k)}\defeq\dfrac{\partial \vecF^{-1, (k)}}{\partial \vecu^{(k)}_{\hg}}=\beta_{\gamma} \dfrac{\vecu^{(k)}_{h,\gamma}}{|\vecu^{(k)}_{h,\gamma}|}.
 \end{align*}\ese
The non-linear variational form~\eqref{discrete_mixed_formulation_interface} may be rewritten in the following canonical form:
find $\vecu_{h,\gamma}\in \mathbf{V}_{\hg}$ and $p_{h,\gamma}\in  M_{\hg}$,
such that
\begin{align*}
 \mathcal{F}_{\gamma}\left[(\vecu_{h,\gamma},p_{h,\gamma}),(\vecv_{\gamma},q_{\gamma})\right]=0, \quad \forall \vecv_{\gamma}\in\mathbf{V}_{h,\gamma},\,\, \forall q_{\gamma}\in M_{h,\gamma},
\end{align*}
where $\mathcal{F}_{\gamma}$ is the residual expression from the mixed system given as follows:
\bse\begin{align*}
&
\mathcal{F}_{\gamma}\left[(\vecu_{h,\gamma},p_{h,\gamma}),(\vecv_{\gamma},q_{\gamma})\right]\defeq
\langle \vecF^{-1}(\vecu_{\hg})\vecu_{h,\gamma}, \vecv_{\gamma} \rangle_{\gamma}+
\langle \mathcal{S}_{\gamma}(p_{\hg}), q_{\gamma}\rangle_{\gamma}\\
&\qquad -\langle p_{\hg}, \nabla_{\tau}\cdot\vecv_{\gamma}\rangle_{\gamma}+
\langle\nabla_{\tau} \cdot\vecu_{h,\gamma}, q_{\gamma}\rangle_{\gamma}
-\langle f_{\gamma}+g_{\gamma}, q_{\gamma}\rangle_{\gamma}.
\end{align*}\ese
In the next step, we calculate the Jacobian given by $\mathcal{J}_{\gamma}\left[(\vecu^{(k)}_{h,\gamma},p^{(k)}_{h,\gamma});(\delta\vecu_{h,\gamma},\delta p_{h,\gamma}),(\vecv_{\gamma},q_{\gamma})\right]$ by taking the G\^{a}teaux
variation of the residual $\mathcal{F}_{\gamma}\left[(\vecu_{h,\gamma},p_{h,\gamma}),(\vecv_{\gamma},q_{\gamma})\right]$ at
$\vecu_{h,\gamma}=\vecu^{(k)}_{h,\gamma}$ and  $p_{h,\gamma}=p^{(k)}_{h,\gamma}$ in the directions of $\delta \vecu_{h,\gamma}$ and $\delta p_{h,\gamma}$, respectively.
This is can be formally obtained by computing
\begin{align*}
&\mathcal{J}_{\gamma}\left[(\vecu^{(k)}_{h,\gamma},p^{(k)}_{h,\gamma});(\delta\vecu_{h,\gamma},\delta
p_{h,\gamma}),(\vecv_{\gamma},q_{\gamma})\right]\\
&\qquad \defeq\Bigg[ \dfrac{\mathcal{F}_{\gamma}\left[(\vecu^{(k)}_{h,\gamma}+\epsilon \delta\vecu_{\hg},p^{(k)}_{h,\gamma}+\epsilon \delta p_{\hg}),(\vecv_{\gamma},q_{\gamma})\right]-\mathcal{F}_{\gamma}\left[(\vecu^{(k)}_{h,\gamma},p^{(k)}_{h,\gamma}),(\vecv_{\gamma},q_{\gamma})\right]}{\epsilon}\Bigg]_{\epsilon\rightarrow 0}.
\end{align*}
This definition yields
\begin{align*}
&\mathcal{J}_{\gamma}\left[(\vecu^{(k)}_{h,\gamma},p^{(k)}_{h,\gamma});(\delta\vecu_{h,\gamma},\delta
 p_{h,\gamma}),(\vecv_{\gamma},q_{\gamma})\right]\\
&\qquad\defeq\langle (\vecF^{-1, (k)} + \vecF_{\partial}^{-1, (k)}\otimes
 \vecu_{\hg}^{(k)})\delta\vecu_{h,\gamma}, \vecv_{\gamma}\rangle_{\gamma}
 +\langle\mathcal{S}_{\gamma}(\delta p_{\hg}), q_{\gamma}\rangle_{\gamma}\\
 &\qquad\quad-\langle\delta p_{\hg}, \nabla_{\tau}\cdot\vecv_{\gamma}\rangle_{\gamma}+
 \langle \nabla_{\tau} \cdot\delta \vecu_{h,\gamma}, q_{\gamma}\rangle_{\gamma},
\end{align*}
 where $\otimes$ denotes the standard tensor product. In each Newton iteration, we solve the following linear variational problem: find
$\delta \vecu_{h,\gamma}\in \mathbf{V}_{h,\gamma}$ and $\delta p_{h,\gamma}\in M_{h,\gamma}$, such that
\begin{align}
\label{discrete_mixed_jacobian_iterative}&\mathcal{J}_{\gamma}\left[(\vecu^{(k)}_{h,\gamma},p^{(k)}_{h,\gamma});(\delta\vecu_{h,\gamma},\delta p_{h,\gamma}),(\vecv_{\gamma},q_{\gamma})\right]=- \mathcal{F}_{\gamma}\left[(\vecu^{(k)}_{h,\gamma},p^{(k)}_{h,\gamma}),(\vecv_{\gamma},q_{\gamma})\right], \quad \forall ( \vecv_{\gamma},q_{\gamma})\in M_{h,\gamma}\times\mathbf{V}_{h,\gamma}.
\end{align}
In compact form, the linear system for the Newton step has the following mixed structure
\begin{align}\label{compact_system_jacobian}
\begin{bmatrix}
 \mathcal{J}_{\gamma}^{k}& B_{\gamma}^{\textn{T}}\\[2mm]
 B_{\gamma}& \mathcal{S}_{\gamma}\\[2mm]
\end{bmatrix}\begin{bmatrix}
 \delta\vecu_{\hg}\\[2mm]
 \delta p_{\hg}\\[2mm]
 \end{bmatrix} =  \begin{bmatrix}
 \mathcal{R}^{\vecu}_{\gamma}\\[2mm]
 \mathcal{R}^{p}_{\gamma}
\end{bmatrix}.
\end{align}

 The interface system~\eqref{compact_system_jacobian} is then solved with  the
 GMRes iterations~\eqref{gmres_pseudocode}. On each GMRes iteration,
 we need to evaluate the action of the Robin-to-Neumann operator
 $\mathcal{S}_{\gamma}$ using Algorithm~\ref{Eval_op}. The solution of the
 interface problem is therefore obtained in an iterative fashion using the following update equations
 until the Newton residual reaches some prescribed tolerance:
 \bse\begin{align*}
&\vecu_{\hg}^{(k+1)}= \vecu_{\hg}^{(k)}+\delta \vecu_{\hg},\\
&p_{\hg}^{(k+1)}= p_{\hg}^{(k)}+\delta p_{\hg}.
\end{align*}\ese
The  result of this iterative approach is then used to infer
the solution in the sub-domains using~\eqref{subdomain_solutions}, which needs
two additional  sub-domain solves.
%
\begin{rem}[An alternative to~\eqref{discrete_mixed_jacobian_iterative}]
For the mixed Jacobian problem in the
fracture~\eqref{discrete_mixed_jacobian_iterative}, it is
possible to adopt the idea introduced in Remark~\ref{remark:FVCGscheme} to
reduce the computational cost by
reformulating~\eqref{discrete_mixed_jacobian_iterative}  into  a cell-centered
finite volume problem with the pressure step $\delta p_{\hg}$ as the sole
variable. The resulting system is also symmetric definite and positive and   can
be solved with   the CG method equipped with a local preconditioner.
\end{rem}

\section{Outer--inner interface iterations with multiscale flux basis}
\label{sec:Mufbi}
As noticed previously,   the dominant computational cost in the above
linearization--domain-decomposition procedures comes from the sub-domain solves   to evaluate the action of $\mathcal{S}_{\gamma}$  using Algorithm~\ref{Eval_op} (step 2(b)). We recall that
the number of sub-domain solves required by each method   is
approximately equal to  $\sum_{k=1}^{N_{\textnormal{lin}}} N^{i}_{\textnormal{dd}}$,
where $N_{\textnormal{lin}}$ is the number of  iterations of the linearization procedure,
and  $N^{k}_{\textnormal{dd}}$ denotes the number of domain decomposition
iterations (GMRes or any Krylov solver).  Even though all sub-domain solves can be
computed in parallel,  this  still be very costly; first,  as the non-linear interface solver may converge very slowly and, second, that  at each linearization iteration the condition number of the linearized interface problem~(\eqref{discrete_mixed_formulation_iterative} for Method~1 and \eqref{discrete_mixed_jacobian_iterative} for Method~2) is large due to a
highly refined mesh.

One way to  reduce the computational costs, is to employ the multiscale  flux
basis, following~\cite{MR2557486}. The motivation of these techniques in this
work stems from eliminating the dependency     between the total number of
solves and the employed  outer-inner procedure  on the
interface-fracture. This is easily achieved by  pre-computing and storing the
flux  sub-domain responses, called \textit{multiscale flux basis}, associated with each
fracture pressure degree of freedom on each sub-domain.

The multiscale flux basis requires solving a fixed
number of linear sub-domain solves and permits retrieving  the action of
$\mathcal{S}_{\gamma}$ on $M_{\hg}$ by simply taking a linear combination of
multiscale flux basis functions. As a result, the number of sub-domains solves
is now independent of the used linearization procedure as well as of the used
solver for the inner domain decomposition systems. In practice,   the number of sub-domains solves
will be reduced if  $\sum_{k=1}^{N_{\textnormal{lin}}} N^{k}_{\textnormal{dd}}$
exceeds the maximum number of fracture pressure degrees of freedom on any sub-domain.

\subsection{Multiscale flux basis }
Following~\cite{MR2557486}, we define
$(\Phi^{\ell}_{\hg})^{\mathcal{N}_{\hg}}_{\ell=1}$ to be the  set of basis
functions on the  interface pressure space $M_{\hg}$, where $\mathcal{N}_{\hg}$
is the number of pressure degrees of freedom on sub-domain $\gamma$. As a result,
on the fracture-interface, we have
\begin{align*}
 p_{\hg}\defeq \sum_{\ell=1}^{\mathcal{N}_{\hg}} p^{\ell}_{\hg}\Phi^{\ell}_{\hg}.
\end{align*}
We compute  the   multiscale flux basis functions corresponding to $(\Phi^{\ell}_{\hg})^{\mathcal{N}_{\hg}}_{\ell=1}$ using the following algorithm:
\begin{algo}[Assembly of the multiscale flux basis]\label{Eval_MS}~
{
\setlist[enumerate]{topsep=0pt,itemsep=-1ex,partopsep=1ex,parsep=1ex,leftmargin=1.5\parindent,font=\upshape}
\begin{enumerate}
    \item Enter the basis $(\Phi^{\ell}_{\hg})^{\mathcal{N}_{\hg}}_{\ell=1}$. Set $\ell=0$.
    \item \textn{\textbf{Do}}
    \begin{enumerate}
        \item Increase $\ell\eq \ell +1$.
        \item Project $\Phi^{\ell}_{\hg}$ on the  sub-domain boundary, i.e,
            $$ \Phi^{\ell}_{\hg} \overset{\mathcal{D}_{h,i}}{\longrightarrow}\lambda_{h,i}^{\ell}.$$
        \item Solve  problem~\eqref{subdomain_problems} in each sub-domain $\Omegai$
            with Robin  boundary condition~$\lambda_{h,i}^{\ell}$ and with $f=0$.
        \item Project the boundary flux onto the mortar space on the
            fracture, i.e,
            $$ -\vecu_{h}(\lambda_{h,i}^{\ell},0)\cdot\vecn_{i} \overset{\mathcal{D}^{\textn{T}}_{h,i}}{\longrightarrow} \Psi_{\hg,i}^{\ell}  $$
    \end{enumerate}
    \textn{\textbf{While} $\ell\leq \mathcal{N}_{\hg}$.}
    \item Form  the  multiscale  flux  basis  for
        sub-domain $\Omegai$, i.e, $$\left\{\Psi_{\hg,i}^{1},\Psi_{\hg,i}^{2},\cdots, \Psi_{\hg,i}^{\mathcal{N}_{\hg}}\right\}\subset M_{\hg}.$$
\end{enumerate}
}
\end{algo}
Once the multiscale flux basis functions are constructed for each sub-domain,
the action of interface operator
$\mathcal{S}^{\textn{RtN}}_{i}$, and then also the action of
$\mathcal{S}_{\gamma}$ via~\eqref{right_hand_side}, is replaced by a linear combination of the multiscale flux basis functions $\Psi_{\hg,i}^{\ell}$. Specifically,
for an interface datum $\varphi_{\hg}\in M_{\hg}$, we have
$\varphi_{\hg}=\sum_{\ell=1}^{\mathcal{N}_{\hg}} \varphi^{\ell}_{\hg}\Phi^{\ell}_{\hg}$, and for $i=1,2,$
\begin{align*}
    \mathcal{S}^{\textn{RtN}}_{i}(\varphi,0)&\eq\mathcal{S}^{\textn{RtN}}_{i}(\sum_{\ell=1}^{\mathcal{N}_{\hg}} \varphi^{\ell}_{\hg}\Phi^{\ell}_{\hg},0)\\
    &\eq\sum_{\ell=1}^{\mathcal{N}_{\hg}} \varphi^{\ell}_{\hg}\mathcal{S}^{\textn{RtN}}_{i}(\Phi^{\ell}_{\hg},0)=\sum_{\ell=1}^{\mathcal{N}_{\hg}} \varphi^{\ell}_{\hg} \Psi^{\ell}_{\hg,i}.
\end{align*}
\begin{rem}
 We observe that  each fracture pressure basis function
 $\Phi^{\ell}_{\hg}$ on the fracture-interface corresponds to exactly two
 different multiscale flux basis functions, one for $\Omega_{1}$ and one for
 $\Omega_{2}$. For the case of a fractures network, say $\gamma=\cup_{i\neq j}
 \gamma_{ij}$, where $\gamma_{ij}$ is the fracture between the sub-domain
 $\Omega_{i}$ and $\Omega_{j}$, the previous basis reconstruction is then
 applied independently on each fracture.
\end{rem}

\subsection{Application on intersecting fractures model: solving the DFNs system}\label{Appendix1}
In this part, we first introduce and describe the case of intersecting fractures, and
then   we provide  our amendments to the previous algorithms.

\subsubsection{Mathematical model}
For the sake of simplicity, we consider the Darcy-Forchheimer model in a two-dimensional geological  domain  made up with three sub-domains $\Omegai$, $i=1,2,3$, physically
subdivided by fractures $\gamma_{i,j}$, $1\leq i<j\leq 3$. The rock matrix is now defined as $\overline{\Omega}=\sum_{i=1}^{3}\overline{\Omega_{i}}$,
$\Omega_{i}\cap \Omega_{j}=\emptyset$, where a single fracture is
$\gamma_{i,j}=\partial \Omega_{i}\cap\partial \Omega_{j}$,  all fractures that
touch sub-domain $\Omegai$ are $\gamma_{i}=\partial \Omega_{i}\setminus \partial
\Omega$. Also, $T=\partial\gamma_{1,2} \cap
\partial\gamma_{2,3}=\partial\gamma_{2,3} \cap
\partial\gamma_{1,3}=\partial\gamma_{1,3} \cap \partial\gamma_{1,2}$ corresponds
to the intersection point of the fractures $ \gamma_{i,j} $ and
$\Gamma_{i}=\partial \Omega_{i}\cap \partial \Omega$ the boundary of each sub-domain $ \Omega_i $.
\begin{figure}[htb]
    \centering
    \resizebox{0.2\textwidth}{!}{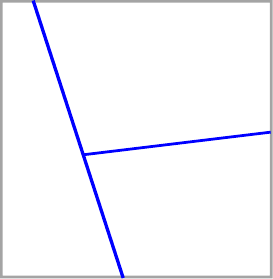}\hspace*{0.05\textwidth}%
    \includegraphics[width=0.2\textwidth]{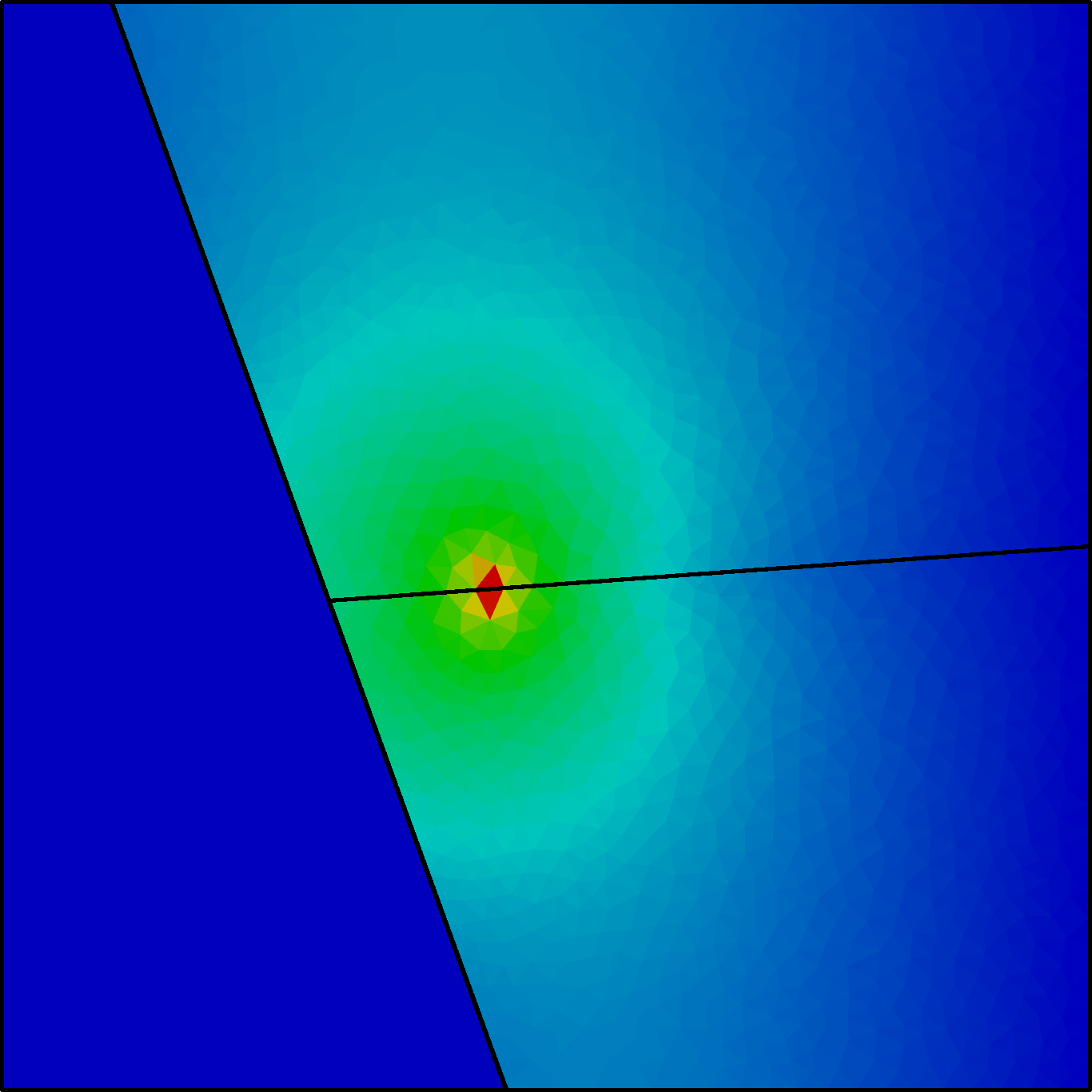}
    \caption{On the left, graphical example of problem
    \eqref{Initial_system_porous}-\eqref{Initial_system_r_interface} along with
    \eqref{eq:intersect_coupling} in case of
    intersecting fractures. On the right, example of construction of a
    multiscale flux basis.}%
    \label{fig:intersecting_domain}
\end{figure}
We impose the Darcy model~\eqref{Initial_system_porous} in each sub-domain $\Omegai$ and the Darcy-Forchheimer
model~\eqref{Initial_system_fracture} in each fracture $\gamma_{i,j}$, with unknowns denoted by
$(\vecu_{\gamma_{i,j}},p_{\gamma_{i,j}})$. See
Figure~\ref{fig:intersecting_domain} (left) as an example.

They are coupled using
the Robin boundary conditions  given  by
\begin{align}\label{Initial_system_r_interface_m}
-\vecu_{i}\cdot\vecn_{i}+\alpha_{i,j} p_{i} =\alpha_{i,j} p_{\gamma_{i,j}}
\quad \textn{on} \;   \gamma_{i,j},
\end{align}
for $1\leq i<j\leq 3$,
where the coefficient  $\alpha_{i,j}$ can now be different in each fracture.
To close the system, we need to impose transmission conditions between the fractures at
the $(d-2)$-dimensional interface $T$. On the intersection $T$, we set
\bse\label{eq:intersect_coupling}\begin{alignat}{2}
    &-\vecu_{\gamma_{i,j}}\cdot\vecn_{{i,j}}+\alpha_{\gamma_{i,j}} p_{\gamma_{i,j}}
    =\alpha_{\gamma_{i,j}} p_{T}\quad&& \textn{on} \; T,\\
    &\sum_{1\leq i<j\leq 3}\vecu_{\gamma_{i,j}}\cdot\vecn_{i,j}=0\quad &&\textn{on} \; T,
\end{alignat}\ese
for $1\leq i<j\leq 3$ where $\vecn_{i,j}$ is the outer unit normal vector to $\partial\gamma_{i,j}.$

For the partition of the sub-domain $\Omegai$, $1\leq i\leq 3$,  and the fractures $\gamma_{i,j}$, $1\leq i<j\leq 3$,
we extend the notation introduced  in
Subsection~\ref{sec:problem_form.subsec:mfe}. We let $\Tau_{\hi}$ be a
partition of the sub-domain
$\Omegai$ into $2$-dimensional simplicial elements and  let  $\Tau_{h, \gamma_{i,j}}$ to be a partition of the fracture $\gamma_{i,j}$ into
$1$-dimensional simplicial elements.   Again, the meshes $\Tau_{\hi}$, $1\leq
i\leq 3$, are  allowed to be non-conforming on the
fractures $\gamma_{i,j}$, $1\leq i<j\leq 3$, but also different
from those used in  $\gamma_{i,j}$, $1\leq i<j\leq 3$ (see Figure~\ref{fig:mortar} for more details).  We also extend the same notation for the
approximation spaces in the sub-domains and in the fractures, and additionally we let   $M_{h, T}$ be the space endowed with  constant functions on $T$.

\subsubsection{Domain decomposition formulation}
The extension of the reduced interface problem~\eqref{discrete_mixed_formulation_interface2} to the
present intersecting fractures setting is as follows:  find the triplet $(\vecu_{\hg},p_{\hg},p_{h,T})\in \mathbf{V}_{\hg}\times M_{\hg} \times M_{h,T}$
such that,   for each $1\leq i<j\leq 3$,
\bse\label{discrete_mixed_formulation_interface_m}\begin{alignat}{2}
\nonumber &\langle
\vecF^{-1}(\vecu_{h,\gamma})\vecu_{h,\gamma},
\vecv_{\gamma} \rangle_{\gamma_{i,j}}+\alpha_{\gamma_{i,j}}^{-1}\langle\vecu_{h,\gamma}\cdot \vecn_{i,j},
\vecv_{\gamma}\cdot \vecn_{i,j} \rangle_{T}-
\langle p_{\hg}, \nabla_{\tau}\cdot\vecv_{\gamma}
\rangle_{\gamma_{i,j}}\\
&\qquad\qquad\qquad\qquad\qquad\qquad\qquad\qquad\,\,\,\,= -\langle p_{h,T}, \vecu_{\gamma_{i,j}}\cdot\vecn_{i,j}\rangle_T\quad &&\forall \vecv_{\gamma} \in
\mathbf{V}_{h, \gamma_{i,j}},\\
&\langle \nabla_{\tau} \cdot\vecu_{h,\gamma}, q_{\gamma}
\rangle_{\gamma_{i,j}}+\langle \mathcal{S}_{\gamma_{i,j}}(p_{\hg}),
q_{\gamma}\rangle_{\gamma_{i,j}} =\langle f_{\gamma_{i,j}} + g_{\gamma_{i,j}},
q_{\gamma}\rangle_{\gamma_{i,j}}\quad&&\forall   q_{\gamma}\in
 M_{h, \gamma_{i,j}},\\
&\sum_{1\leq i<j\leq 3}\langle \vecu_{\gamma_{i,j}}\cdot\vecn_{i,j}, q_{T}
\rangle_T=0\quad&& \forall q_{T}\in M_{h, T}.
\end{alignat}\ese
 On each fracture, the Robin-to-Neumann operator $\mathcal{S}_{\gamma_{i,j}}$ and the linear functional $g_{\gamma_{i,j}}$, $1\leq i<j\leq 3$, are now given by
 \bse\begin{align*}
  &\mathcal{S}_{\gamma_{i,j}}(p_{\hg})\eq \sum_{l\in(i,j)}\mathcal{S}^{\textn{RtN}}_{\gamma_{l}}(p_{\hg},0)= -\sum_{l\in(i,j)}\vecu_{l}(p_{\gamma},0)\cdot\vecn_{l}|_{\gamma_{l}},\\
  & g_{\gamma_{i,j}}\eq\sum_{l\in(i,j)}\mathcal{S}^{\textn{RtN}}_{\gamma_{l}}(0,f_{l})= \sum_{l\in(i,j)}\vecu_{l}(0,f_{l})\cdot\vecn_{l}|_{\gamma_{i}}.
 \end{align*}\ese
 The above  problem can be seen as a  DFNs system on the set of fractures, and as a domain decomposition problem between the $1$-dimensional fractures $\gamma_{i,j}$, $1\leq i<j\leq 3$, cf.\cite{MR2773340,MR2890281,MR3757111} for more details.

 \subsubsection{Iterative procedure}
 We propose  to solve the  non-linear domain decomposition problem~\eqref{discrete_mixed_formulation_interface_m} using the fixed-point approach in
 Subsection~\ref{subsec:fixed_point}. This iterative process is now equipped with the  multiscale flux  basis of Section~\ref{sec:Mufbi} to lessen the interface iterations. To this aim, we  introduce
  \begin{align*}
  &\mathcal{S}_{\gamma}(p_{\hg})\eq\sum_{0\leq i< j\leq 3}\mathcal{S}_{\gamma_{i,j}}(p_{\hg})\quad\text{and}\quad g_{\gamma}\eq\sum_{1\leq i<j\leq 3}g_{\gamma_{i,j}},
 \end{align*}
 and let
 \begin{align*}
  &\mathcal{S}_{T}(p_{h,T})\eq\sum_{1\leq i<j\leq 3}\vecu_{\gamma_{i,j}}\cdot\vecn_{i,j}|_{T}.
 \end{align*}

Applying the fixed-point algorithm on the set of interface Darcy-Forchheimer  equations~\eqref{discrete_mixed_formulation_interface_m} can be interpreted as follows:
 at the iteration $k\geq 1$, we solve
\begin{align}\label{linearized_dfns_system}
\begin{bmatrix}
 \vecF^{-1, (k)}_\gamma& B_\gamma^{\textn{T}}&\mathcal{S}^{\textn{T}}_{T}\\[2mm]
 B_\gamma& \mathcal{S}_{\gamma}&0\\[2mm]
 \mathcal{S}_{T}& 0&0\\[2mm]
\end{bmatrix}
\begin{bmatrix}
 \vecu^{k}_{\hg}\\[2mm]
 p^{k}_{\hg}\\[2mm]
 p^{k}_{h,T}
 \end{bmatrix} =  \begin{bmatrix}
 0\\[2mm]
 f_{\gamma}+g_{\gamma}\\[2mm]
 0
\end{bmatrix},
\end{align}
using  GMRes method until a fixed tolerance is reached.  Again, the evaluation of $\mathcal{S}_{\gamma}$ in each interface GMRes iteration dominates the total computational costs of this outer--inner procedure. Note that each inner iteration also requires the evaluation of the Dirichlet-to-Neumann operator $\mathcal{S}_{T}$, which requires solves in the fractures. The complete algorithm when equipped with multiscale flux basis is now given by the following algorithm.

\begin{algo}[Fixed-point algorithm with multiscale flux basis for  fracture network  model]\label{Eval_MS}~
{
\setlist[enumerate]{topsep=0pt,itemsep=-1ex,partopsep=1ex,parsep=1ex,leftmargin=1.5\parindent,font=\upshape}
\begin{enumerate}
    \item Enter the source terms and  the permeabilities in the fractures and the rock matrices.
    \item Choose the meshes  $\Tau_{\hi}$, $1\leq i\leq 3$, and $\Tau_{h, \gamma_{i,j}}$, $1\leq i<j\leq 3$.
    \item Calculate the right-hand-sides $g_{\gamma_{i,j}}$, $1\leq i<j\leq 3$, by solving
        the Darcy sub-domain problem in $\Omegai$ with source term $f_{i}$ and
        zero Robin value on the fracture-interface $\gamma_{i}$. Then, compute
        the resulting jump across all sub-domain interfaces.
    \item In the sub-domain $\Omegai$, $1\leq i\leq 3$, let $\mathcal{N}_{h,\gamma_{i}}$ be  the number of
          degrees of freedom in the space $M_{h,\gamma_{i}}$. Define the
          basis
          $(\Phi^{\ell}_{h,\gamma_{i}})^{\mathcal{N}_{h,\gamma_{i}}}_{\ell=1}$.
          Set $i=0$.\\
    \textn{\textbf{Do}} \{\textn{Assembly of the multiscale flux basis}\}
        \begin{enumerate}
            \item Increase $i\eq i +1$.
            \item Compute  the   multiscale flux basis functions $(\Psi^{\ell}_{h,\gamma_{i}})^{\mathcal{N}_{h,\gamma_{i}}}_{\ell=1}$
            corresponding to $(\Phi^{\ell}_{h,\gamma_{i}})^{\mathcal{N}_{h,\gamma_{i}}}_{\ell=1}$ using
            Algorithm~\ref{Eval_MS}, i.e., $$\Psi^{\ell}_{h,\gamma_{i}}\eq
            \mathcal{S}^{\textn{RtN}}_{\gamma_{l}}(\Phi^{\ell}_{h,\gamma_{i}},0),\quad \ell=1,\cdots, \mathcal{N}_{h,\gamma_{i}}.$$
    \end{enumerate}
    \textn{\textbf{While} $i\leq 3$.}
    \item Given an initial guess $\vecu^{(0)}_{h,\gamma_{i,j}}$, $1\leq i<j\leq
    3$. Set $k=0$.\\
    \textn{\textbf{Do}} \{\textn{Fixed-point iterations}\}
    \begin{enumerate}
        \item Increase $k\eq k +1$.
        \item Solve the linear system on the fractures~\eqref{linearized_dfns_system} using GMRes method~\eqref{gmres_pseudocode},
         where in every iteration the operator action $\mathcal{S}_{\gamma}$  on any $\varphi_{\hg}\in M_{\hg}$ is computed with the
            following steps:
        \begin{enumerate}
            \item  Use a linear combination of the multiscale flux basis to compute the   action of $\mathcal{S}^{\textn{RtN}}_{\gamma_{i}}$ by
                $$\mathcal{S}^{\textn{RtN}}_{\gamma_{i}}(\varphi_{\hg,i},0)\eq \sum_{\ell=1}^{\mathcal{N}_{\hg,i}} \varphi^{\ell}_{\hg,i} \Psi^{\ell}_{\hg,i}.$$
            \item Compute the jump across all the fractures:
                $$\mathcal{S}_{\gamma}(\varphi_{\hg})\eq\sum_{0\leq i< j\leq 3}\sum_{l\in(i,j)}\mathcal{S}^{\textn{RtN}}_{\gamma_{l}}(\varphi_{\hg},0).$$
        \end{enumerate}
    \end{enumerate}
     \textn{\textbf{While}} $\dfrac{\|(p^{k,\infty}_{\hg},\vecu^{k,\infty}_{\hg}) - (p^{k-1,\infty}_{\hg},\vecu^{k-1,\infty}_{\hg})\|_{\infty} }{\|(p^{k-1,\infty}_{\hg},\vecu^{k-1,\infty}_{\hg})\|_{\infty}} \geq \varepsilon_{tol}$.\hfill\eqnum\label{class_stp_critera1}
\end{enumerate}
}
\end{algo}

\section{Numerical examples}
\label{sec:examples}

In this section, we validate the model and analysis presented in the previous
parts by means of numerical test cases. We have chosen three examples designed to show how the proposed linearization--domain-decomposition approaches equipped with    multiscale flux basis behaves vs the standard ones in various physical and geometrical situations. To compare these  approaches, the main criteria considers the number of solutions of the higher-dimensional sub-problems since it constitutes the major computational
cost. We consider solving the problem in the network of fractures as
negligible. Since each of the higher-dimensional sub-problem is linear and  will be solved many times, we consider an LU-factorization of the system matrix and a forward-backward substitution algorithm to compute the numerical solution.
It results in a computational cost reduced to $\mathcal{O}(n^2)$ flops
each time, where $n$ is the size of the matrix. For bigger systems, an
iterative scheme is preferable.

We use the PorePy \cite{Keilegavlen2017a} library, which is a simulation tool
for fractured and deformable porous media written in Python. PorePy uses SciPy
\cite{Jones2001} as default sparse linear algebra.  All the examples are
reported in the GitHub repository of PorePy, we want to stress again that even
if we focus on lowest-order Raviart-Thomas-N\'{e}d\'{e}lec finite elements,
our implementation is agnostic with respect to the numerical scheme.

For the multiscale flux basis scheme presented in Section \ref{sec:Mufbi}, for a fixed rock
matrix grid and normal fracture permeability it is possible to compute once all
the basis functions. The results in the next parts should be read under this
important property of the method, thus in many cases only a pure fracture
network will be solved at a negligible computational cost. The multiscale basis
functions are computed and stored in an \textit{offline} phase prior the
simulation (called \textit{online}).

Unless otherwise noted, the tolerance for the relative residual in
the inner GMRes algorithm   is taken to be   $10^{-6}$. The same tolerance is chosen for the outer Newton/fixed-point algorithm. We consider an LU-factorization of the fracture network matrix~\cite{MR1941848,MR1323819} as the preconditioner of the GMRes method.

In the examples, we use the abbreviation MS when the linearization--domain--decomposition approach is equipped with multiscale flux basis techniques, and DD for the corresponding classical approach.

\begin{rem}[Fracture aperture]
Even if not explicitly considered in the previous parts of the work, we
introduce the fracture aperture $\epsilon$ as a constant parameter. This choice
is based on the fact that geometries and (some) data of the forthcoming examples
are taken from the literature.
\end{rem}

In \ref{subsec:setting} we describe the geometry and some data of the problem
considered. Few subsections follows with an increase level of challenge: linear
case in \ref{subsec:example1}, Forchheimer model in \ref{subsec:example2},
Forchheimer model with heterogeneous parameters in \ref{subsec:example3}, to
conclude with a generalized Forchheimer model in \ref{subsec:example4}.

\subsection{Problem setting}\label{subsec:setting}

\begin{figure}[hbt]
    \centering
    \includegraphics[width=0.35\textwidth]{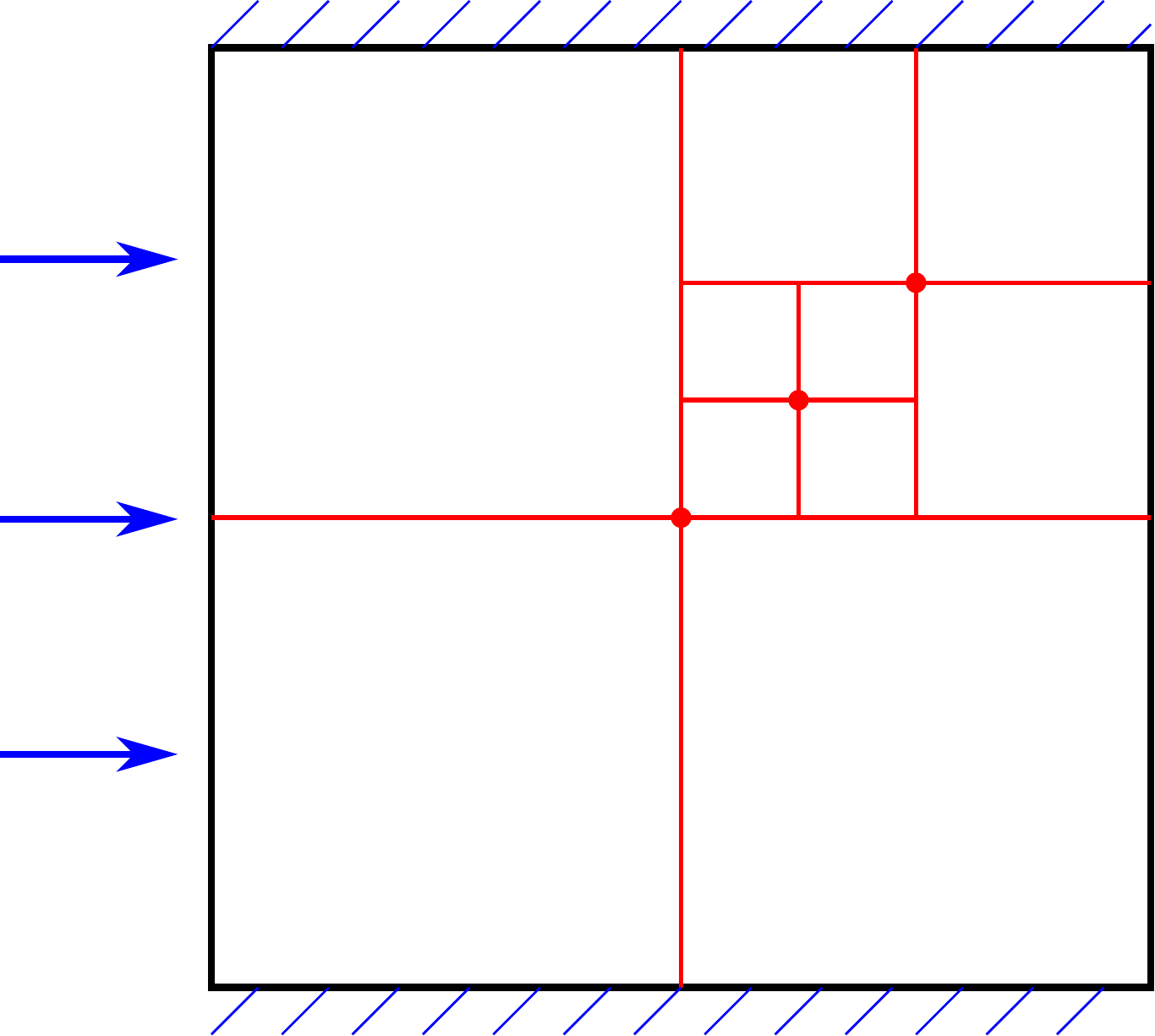}
    \caption{Graphical representation of the domain and fracture network
    geometry common for all test cases.}
\label{fig:geiger_domain}
\end{figure}

To validate the performance of the two proposed algorithms, we consider the
first problem presented in the benchmark study \cite{Flemisch2016a}. The unit
square domain $\Omega$, depicted in Figure \ref{fig:geiger_domain}, has unitary
permeability of the rock matrix and it is divided into  10 sub-domains by a set of fractures
with fixed aperture $\epsilon$ equal to $10^{-4}$.
At the boundary, we impose zero flux condition on the top and bottom, unitary
pressure on the right, and flux equal to $-1$ on the left. The boundary conditions are applied
to both the rock matrix and the fracture network.

Contrary to what has been done in the
benchmark paper, we consider four different scenarios for the fracture
permeabilities, by having high or low values in the tangential and normal parts.
Thus, we have the \textit{case (i)} with high permeable fractures,
\textit{case (ii)} has low permeable fractures, while \textit{cases (iii)} and
\textit{(iv)} have mixed high and low permeability in normal and tangential
directions. See Table~\ref{tab:case1} for a summary of the fracture
permeability in each case.
\begin{table}[htbp]
    \centering
    \begin{tabular}{l|p{1.2cm}|p{1.2cm}|}
        \cline{2-3}
        & \multicolumn{1}{c|}{$\mathbf{K}_\gamma$} &
        \multicolumn{1}{c|}{$\alpha_\gamma$} \\ \hline
        \multicolumn{1}{|l||}{\textit{case (i)}}   & $10^4\epsilon \mathbf{I}$    & $10^4/\epsilon$   \\ \hline
        \multicolumn{1}{|l||}{\textit{case (ii)}}  & $10^{-4}\epsilon\mathbf{I}$ & $10^{-4}/\epsilon$\\ \hline
        \multicolumn{1}{|l||}{\textit{case (iii)}} & $10^4\epsilon\mathbf{I}$    & $10^{-4}/\epsilon$\\ \hline
        \multicolumn{1}{|l||}{\textit{case (iv)}}  & $10^{-4}\epsilon\mathbf{I}$ & $10^4/\epsilon$   \\ \hline
    \end{tabular}
    \caption{Definition of the cases for the examples.}%
    \label{tab:case1}
\end{table}
\textit{Case (i)} and \textit{(ii)} have the same
permeabilities used in the benchmark paper \cite{Flemisch2016a}.

In the following examples, we consider the maximal number of rock matrix solves to be
$10^4$, and we mark with $\infty$ if this is exceeded.

\subsection{Darcy model: $\beta_\gamma = 0$} \label{subsec:example1}

The first example considers the Forchheimer coefficient set to zero, thus the
problem becoming linear. The results for different level of discretization are
reported in Table~\ref{tab:results_case1}. We indicate by \textit{level 1} a
grid with 110 triangles and 26 segments, \textit{level 2} with 1544 triangles
and 84 segments, and \textit{level 3} with 3906 triangles and 138 segments.
\begin{table}[htbp]
    \centering
    \begin{tabular}{l|p{.5cm}|p{.5cm}||p{.5cm}|p{.5cm}||p{.5cm}|p{.5cm}|}
        \cline{2-7}
        & \multicolumn{2}{c||}{\textit{level 1}}
        & \multicolumn{2}{c||}{\textit{level 2}} & \multicolumn{2}{c|}{\textit{level 3}} \\
        \cline{2-7}
        &           MS & DD & MS & DD & MS & DD \\ \hline
        \multicolumn{1}{|l||}{\textit{case (i)}}   & 28$^\dagger$ & 10 & 86$^\dagger$ & 11  & 140$^\dagger$ & 11  \\ \hline
        \multicolumn{1}{|l||}{\textit{case (ii)}}  & 28$^{\mathsection}$ & 81 & 86$^{\mathsection}$ & 112 & 140$^{\mathsection}$ & 189 \\ \hline
        \multicolumn{1}{|l||}{\textit{case (iii)}} & 28$^{\mathsection}$ & 22 & 86$^{\mathsection}$ & 28 & 140$^{\mathsection}$ & 29  \\ \hline
        \multicolumn{1}{|l||}{\textit{case (iv)}}  & 28$^\dagger$ & 82 & 86$^\dagger$ & 61 & 140$^\dagger$ & 86  \\ \hline
    \end{tabular}
    \caption{Total number of  the  higher-dimensional problem solves for the case study of
    example in Subsection \ref{subsec:example1}. For each level of refinement
    cases marked in $\dagger$ share the same multiscale flux basis, which can be
    constructed only once. The same is valid for ${\mathsection}$.}%
    \label{tab:results_case1}
\end{table}

\begin{figure*}[htbp]
    \centering
    \includegraphics[width=0.475\textwidth]{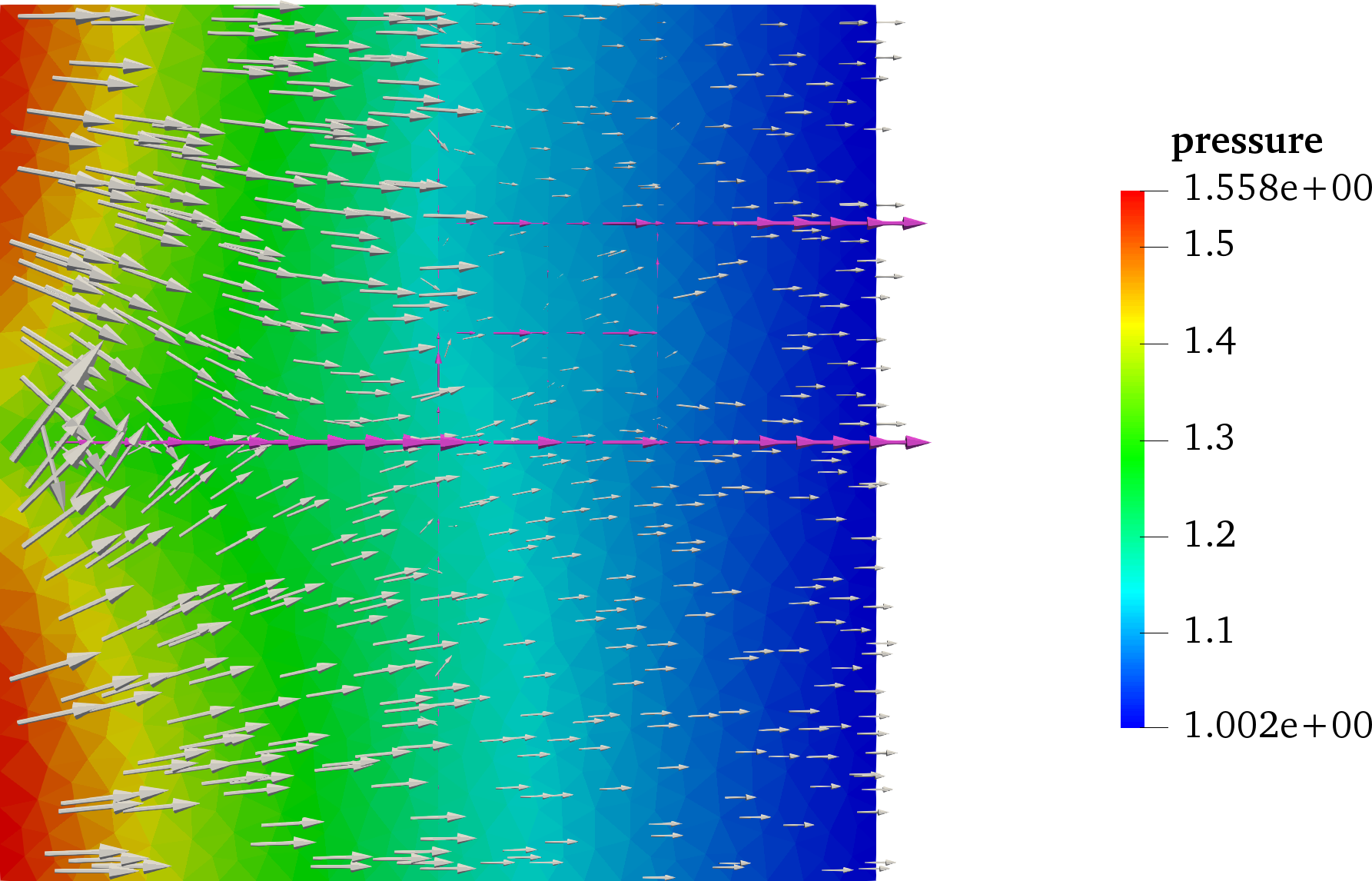}%
    \includegraphics[width=0.475\textwidth]{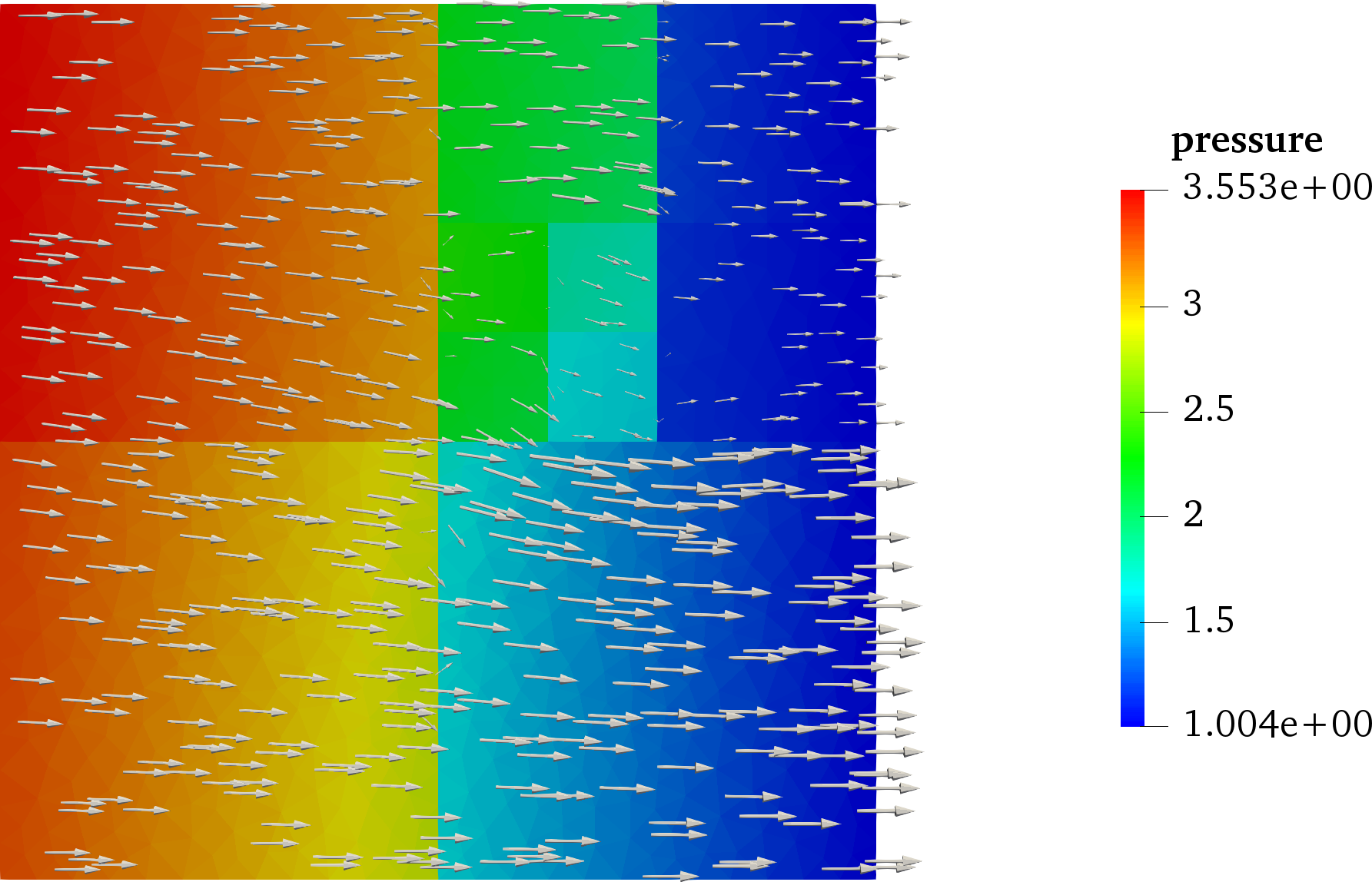}\\%
    \includegraphics[width=0.475\textwidth]{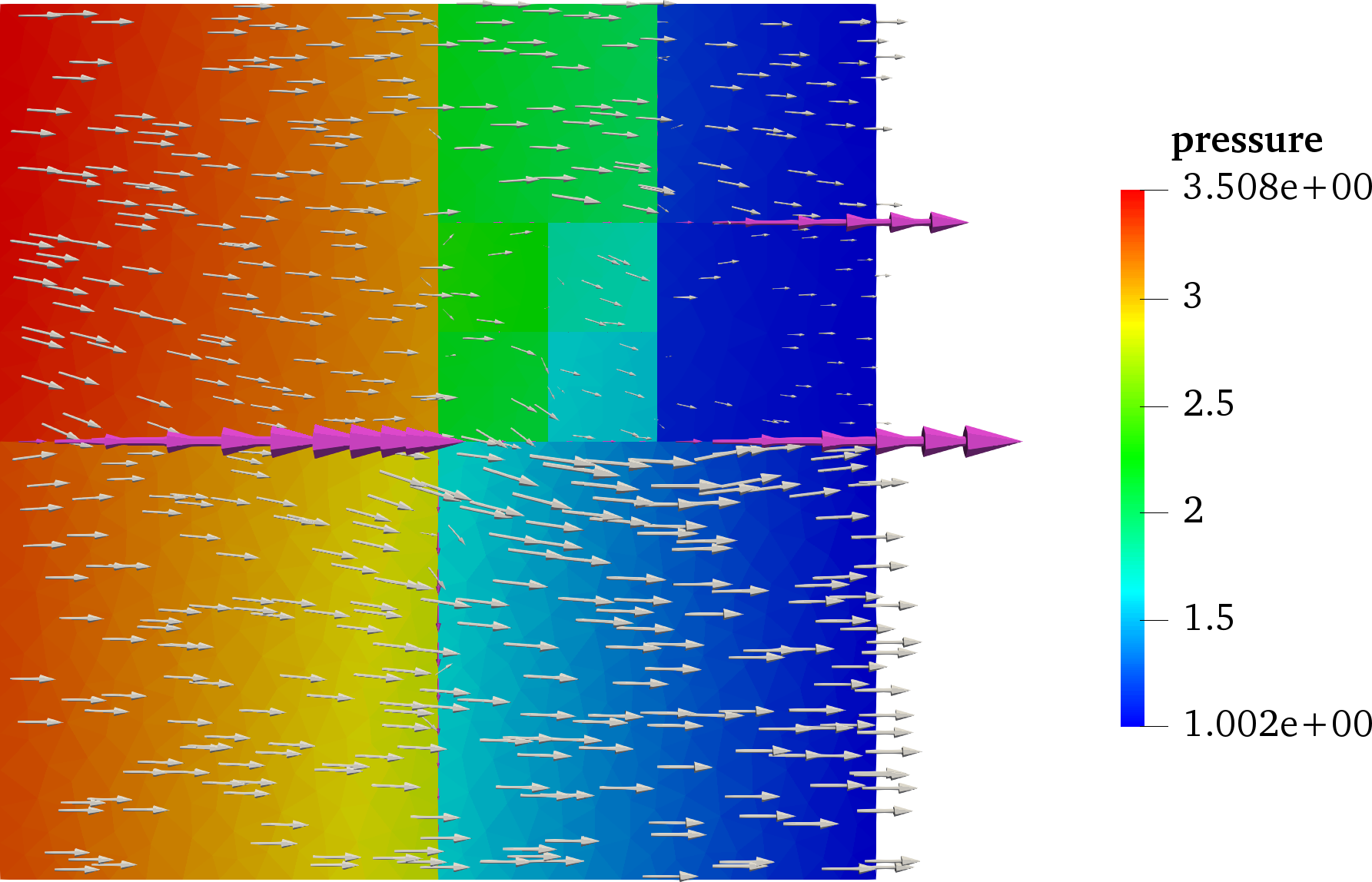}%
    \includegraphics[width=0.475\textwidth]{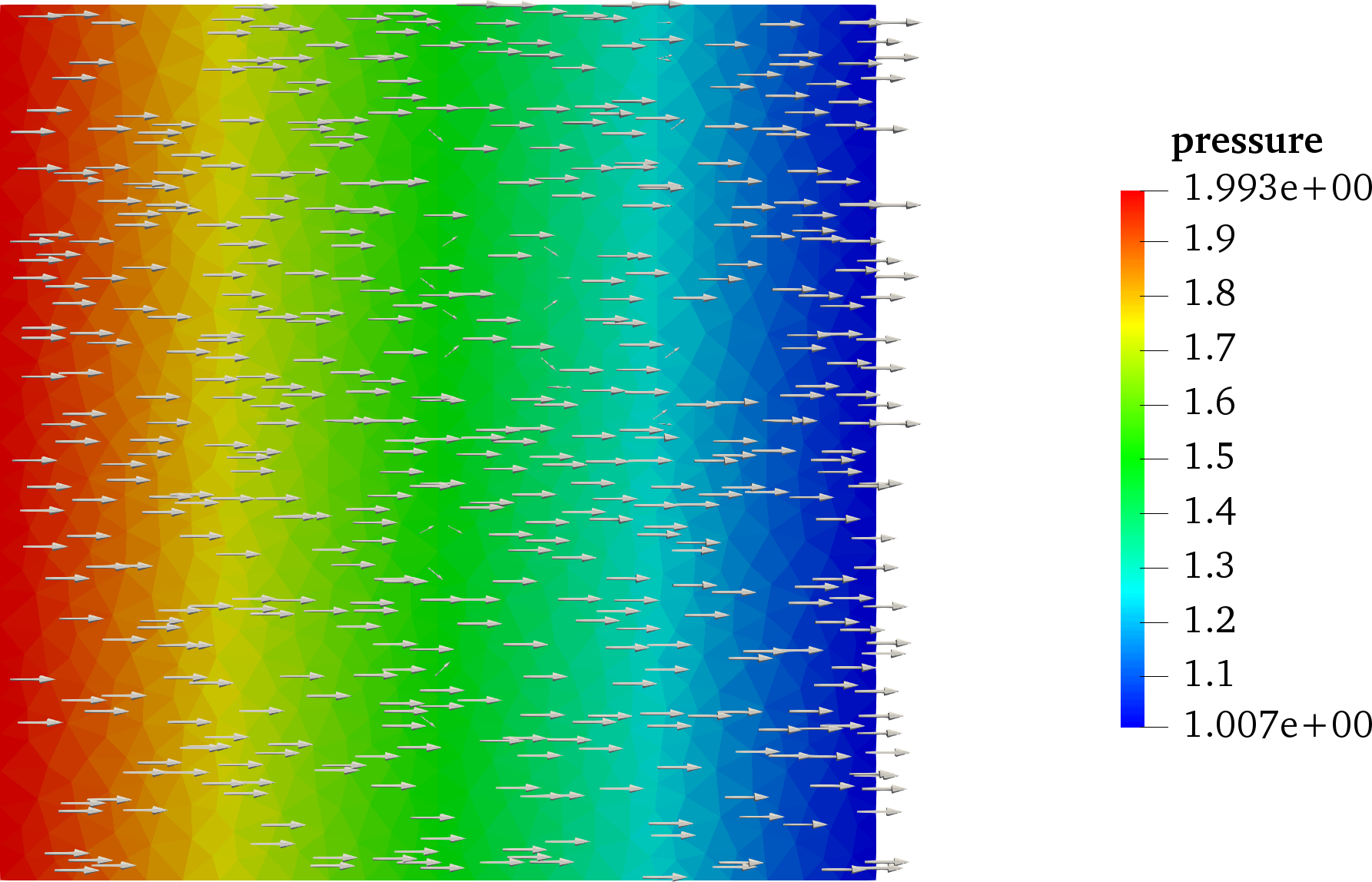}%
    \caption{Pressure and velocity solutions for the four cases: on the top-left
    \textit{case (i)}, on the top-right \textit{case (ii)}, on the bottom-left
    \textit{case (iii)}, and on the bottom-right \textit{case (iv)}.
    In all the
    cases, the velocity is represented by arrows (purple for the fractures)
    proportional to its magnitude.
    }%
    \label{fig:results_case1}
\end{figure*}

Table~\ref{tab:results_case1} shows the results of this example  for the
physical considerations of Table~\ref{tab:case1}. We notice that for high
permeable fractures (\textit{case (i)} and \textit{(iii)}), the standard domain
decomposition method performs better than our method with multiscale flux basis,
while the opposite occurs for low permeable fractures. A possible
explanation is related to the ratio between normal and tangential permeability.
The normal permeability determines how strong the flux exchange is between the
rock matrix and the fractures (thus, the communications at each DD iteration),
while for small values of the tangential permeability the fractures are more
influenced by the surrounding rock matrices. The opposite occurs in the case of
high tangential permeability. Additionally, the choice of the preconditioner for
DD slightly goes in favor of high permeable fractures due to the dominating role
of the fracture flow in the system. We also recall that  the number of
higher-dimensional problem solves does not depend on the number of outer--inner
interface iterations, but only on the  number  of  local  mortar  degrees  of
freedom on the fractures network.  A further important result in this
experiments, is that \textit{case (i)} and \textit{(iv)} share the same value of
$\alpha_\gamma$, thus the multiscale flux basis are computed only once per level
of refinement. The same applies to \textit{case (ii)} and \textit{(iii)}.  As a
result,  the developed method is globally more efficient than the classical
approach. That is, the results in Table~\ref{tab:results_case1} shows a
reduction of the number of the higher-dimensional problem solves from 195  to 56
for \textit{level 1}, from 212 to 186 for \textit{level 2}, and from 312  to 280
for \textit{level 3}.  Note that the two   methods produce the same solution for
all the cases, within the same relative convergence tolerance. The numerical
solution for all cases is reported in Figure \ref{fig:results_case1}.

The next series of numerical experiments aims at assessing the stability of the
domain decomposition approach with respect to   GMRes tolerance. The  multiscale
flux basis approach provides the extra flexibility to do such analysis with
negligible costs, by  reusing the stored multiscale flux basis used for the
results of  Table~\ref{tab:results_case1}  but now with different tolerance for
GMRes. Further,  this set of test cases aims assessing how  the overall gain for an
entire simulation in terms of number of higher-dimensional problem solves  can
be appreciated or depreciated with more  or less stringent stopping criteria for
GMRes; this is a preparatory step to address  the complete
approaches of Section~\ref{sec:nonlinear_algorithms} for the full  nonlinear
problem,  which requires several solves of linear Darcy problems, for which one
should formulate the stopping criteria very carefully. In
Table~\ref{tab:results_case1}, we have considered the relative  residual to be
below $10^{-6}$, while in Table~\ref{tab:results_case1_tol} we present the
results in the case of $10^{-4}$ and $10^{-8}$. Based on the results
of~Table~\ref{tab:results_case1_tol}, we can conclude that even with  less
stringent criterion,  a considerable gain in terms of number of
higher-dimensional problem solves can be achieved. We also see that all the
results are free of oscillations and neither the fracture, barrier, or the very
different  tangential and normal permeabilities  pose any problems for the
domain decomposition approach. Based on the above results  and in what follows
we consider $10^{-6}$ as tolerance for the GMRes algorithm.
{
\renewcommand{\arraystretch}{1.1}
\begin{table}[htbp]
    \centering
    \begin{tabular}{l|p{.5cm}|p{.5cm}||p{.5cm}|p{.5cm}||p{.5cm}|p{.5cm}|}
        \cline{2-7}
        & \multicolumn{2}{c||}{\textit{level 1}}
        & \multicolumn{2}{c||}{\textit{level 2}} & \multicolumn{2}{c|}{\textit{level 3}} \\
        \hline
        \multicolumn{1}{|l||}{tolerance} &  $10^{-4}$ & $10^{-8}$ & $10^{-4}$ &
        $10^{-8}$ & $10^{-4}$ & $10^{-8} $\\ \hline
        \multicolumn{1}{|l||}{\textit{case (i)}}   & 8  & 11  & 9   & 12 & 9 & 12 \\ \hline
        \multicolumn{1}{|l||}{\textit{case (ii)}}  & 42 & 105 & 82  & $\infty$ & 150 & $\infty$ \\ \hline
        \multicolumn{1}{|l||}{\textit{case (iii)}} & 21 & 30  & 22 & 36  & 22 & 36\\ \hline
        \multicolumn{1}{|l||}{\textit{case (iv)}}  & 42 & 122 & 50  & $\infty$ & 70 & $\infty$ \\ \hline
    \end{tabular}
    \caption{Total number of  the higher-dimensional problem solves for the case study of
    example in Subsection \ref{subsec:example1}. For each level of refinement
    we change the convergence tolerance for the domain decomposition method.}%
    \label{tab:results_case1_tol}
\end{table}
}

\subsection{Forchheimer model} \label{subsec:example2}

In this second example we consider \textit{case (i)} and \textit{(iii)} for the
fracture permeabilities since the Forchheimer model requires high permeable
fractures. In this problem, we fix the computational grid \textit{level 2} of
Table~\ref{tab:results_case1} and we change the value of $\beta_\gamma$ in order
to compare the performances of Method 1 and Method 2  with and without multiscale flux basis. The Forchheimer coefficient here varies as $\{1, 10^2, 10^4, 10^6\}$. These values are reasonable since in our model we do not explicitly
scale $\beta_\gamma$ by the aperture, as done in \cite{Frih2008,Knabner2014}. Therefore, the last two values are more realistic. The stopping criteria for both methods is based on the relative residual criteria~\eqref{class_stp_critera1} with a threshold fixed as $10^{-6}$. The initial guess is
taken by  solving the linear Darcy by taking $\beta_\gamma$ equal to zero.
\begin{table}[htbp]
    \centering
    \begin{tabular}{l|l|l||l|l|}
        \cline{2-5}
        & \multicolumn{2}{c||}{\textit{case (i)}} & \multicolumn{2}{c|}{
        \textit{case (iii)}} \\
        \hline
        \multicolumn{1}{|l||}{$\beta_\gamma$} & MS      & DD       & MS & DD    \\ \hline
        \multicolumn{1}{|l||}{$1$}           & 86$^\dagger$ (2)  & 33   (2)  &   86$^{\mathsection}$ (1) & 56 (1)  \\ \hline
        \multicolumn{1}{|l||}{$10^2$}            & 86$^\dagger$ (3)  & 44   (3)        &86$^{\mathsection}$ (2) & 84 (2) \\ \hline
        \multicolumn{1}{|l||}{$10^4$}          & 86$^\dagger$ (8) & 99 (8) &        86$^{\mathsection}$ (3) & 115 (3) \\ \hline
        \multicolumn{1}{|l||}{$10^6$}          & 86$^\dagger$ (94) & 1424 (94) & 86$^{\mathsection}$ (11) &  457 (11) \\ \hline
    \end{tabular}
    \caption{Total number of  the higher-dimensional problem solves required by Method~1 for the case
    study in Subsection \ref{subsec:example2}.  The number   of the fixed-point  iterations are in brackets. Within each case
    the construction of the multiscale flux basis is done only once, we mark by
    $\dagger$ (respectively ${\mathsection}$) common computations.}%
    \label{tab:results_case2}
\end{table}

\textbf{For Method 1}, the number of   higher-dimensional problem solves is reported in Table~\ref{tab:results_case2}. As expected, Method~1 equipped with multiscale flux basis (MS)  performs all the higher-dimensional problem solves  in the offline phase, thus the outer--inner interface iterations for the resulting fracture network problem do not influence the total computational costs. On the
contrary, the computational costs of the
classical approach (DD)  is influenced by the non-linearity, by varying
$\beta_{\gamma}$,  as well as by the ratio of the normal and tangential permeabilities, by varying $\mathbf{K}_\gamma$ and $\alpha_\gamma$. Particularly,  the total gain of the new approach is more significant when  the non-linear effects becomes more important (by increasing the value of $\beta_\gamma$).  Furthermore, for the entire simulation of  each case of Table~\ref{tab:results_case2},   the multiscale flux basis  are  computed only once. As a conclusion, the entire simulation  of \textit{case (i)}  required  for Method~1  1600 higher-dimensional problem solves, while for Method~1   with multiscale flux basis,  this number is reduced   by $95\%$. For \textit{case (iii)}, we reduce the computational costs by $88\%$.



\begin{figure*}[tbp]
    \centering
    \includegraphics[width=0.475\textwidth]{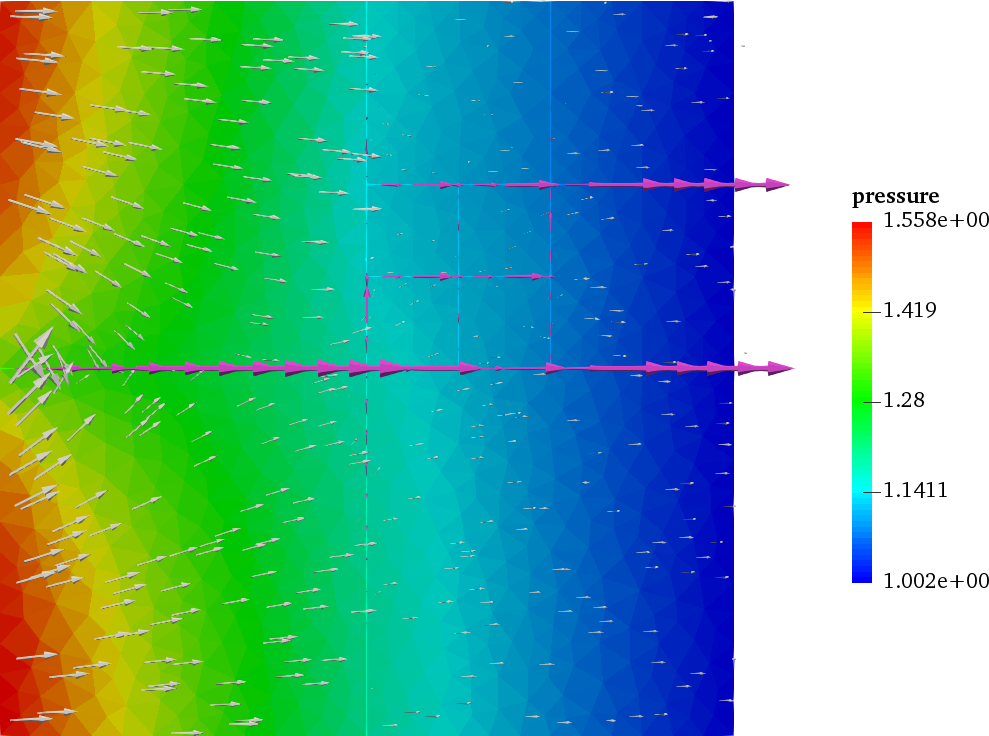}%
    \includegraphics[width=0.475\textwidth]{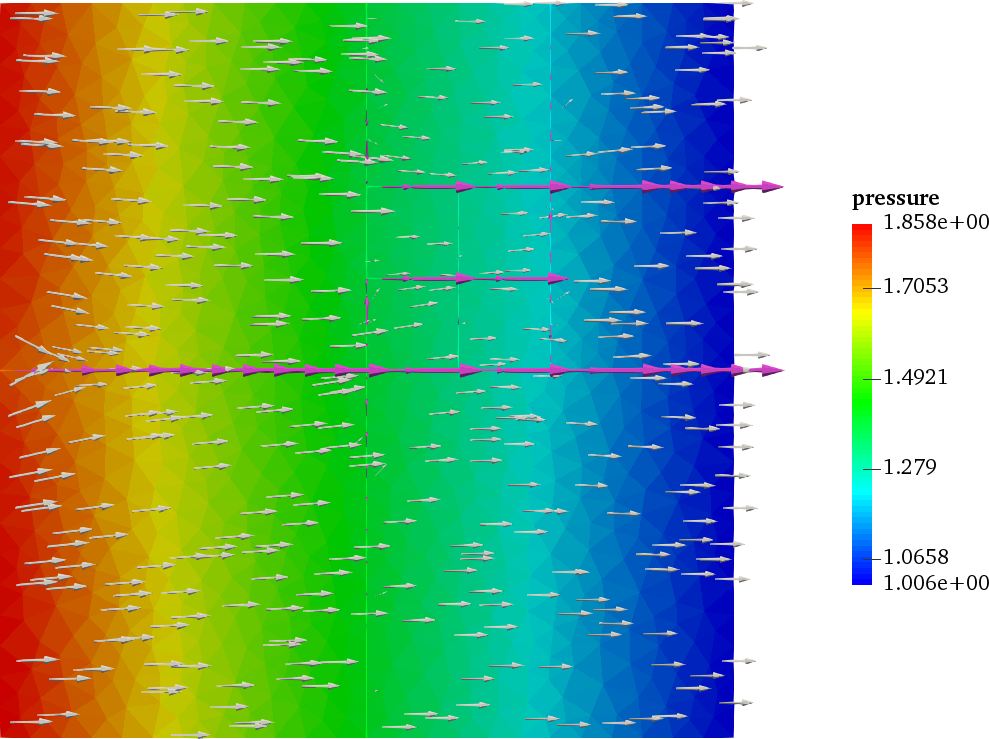}\\%
    \includegraphics[width=0.475\textwidth]{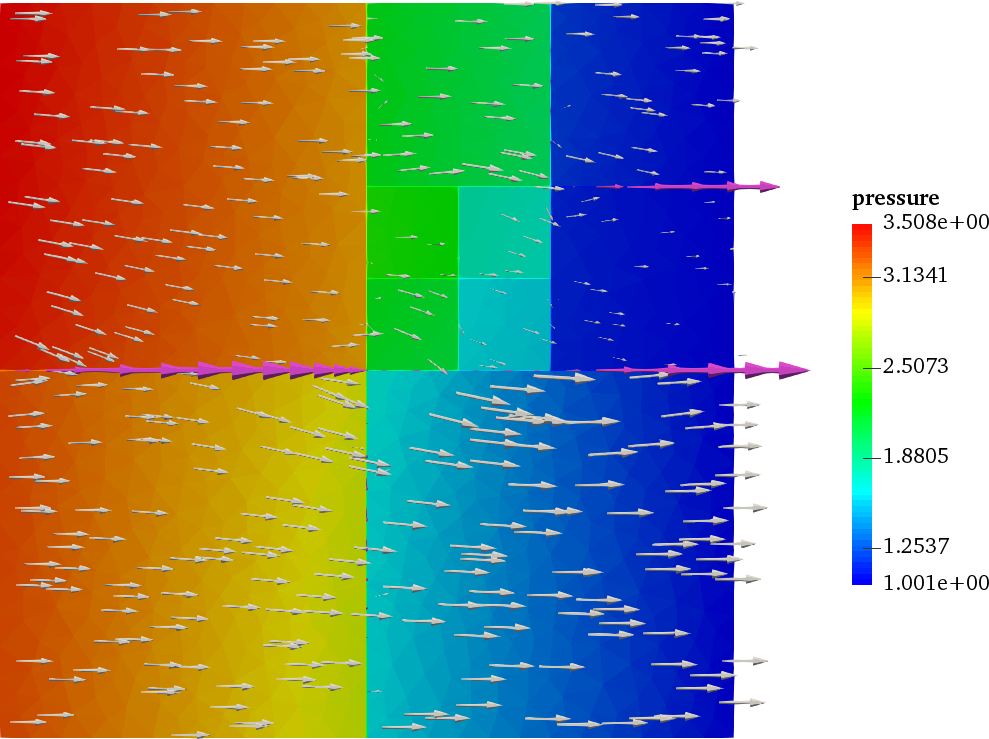}%
    \includegraphics[width=0.475\textwidth]{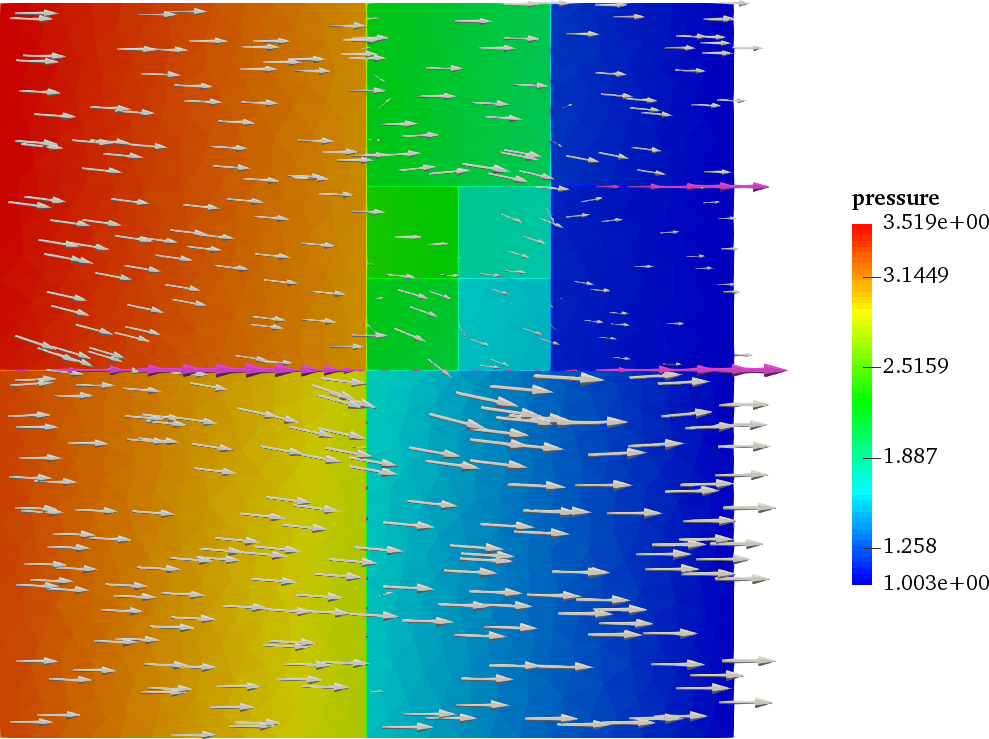}
    \caption{Pressure and velocity solutions for different configurations of
    example presented in Subsection \ref{subsec:example2}. On the top
    \textit{case (i)} and on the bottom for \textit{case (iii)}. On the left, we consider value of the
    Forchheimer coefficient equal to $\beta_\gamma = 1$ and on the right a high
    value $\beta_\gamma = 10^6$. In all the
    cases, the velocity is represented by arrows (purple for the fractures)
    proportional to its magnitude.
    }%
    \label{fig:results_case2}
\end{figure*}
The numerical solution for two values of $\beta_\gamma$ are reported in Figure
\ref{fig:results_case2} for both cases. Despite the different values of
$\beta_\gamma$, we notice that the graphical results are
very similar in the case of low $\alpha_\gamma$. While for high value of
$\alpha_\gamma$, the resulting apparent permeability given by
$\mathbf{K}_\gamma
(1+\mathbf{K}_\gamma^{-1} \beta_\gamma |\mathbf{u}_\gamma|)^{-1}$ decreases (for
a fixed $|\mathbf{u}_\gamma|$) and the fractures are less prone to be the main
path for the flow. Also as stated previously,  since we do not explicitly scale $\beta_\gamma$ by the aperture, values of $\beta_\gamma > 10^{4}$ are more likely
for real applications.

\begin{table}[htbp]
    \centering
    \begin{tabular}{l|l|l||l|l|}
        \cline{2-5}
        & \multicolumn{2}{c||}{\textit{case (i)}} & \multicolumn{2}{c|}{
            \textit{case (iii)}} \\
        \hline
        \multicolumn{1}{|l||}{$\beta_\gamma$} & MS      & DD       & MS & DD    \\ \hline
        \multicolumn{1}{|l||}{$1$}           & 86$^\dagger$ (2)  & 20   (2)  &   86$^{\mathsection}$ (1) & 38 (1)  \\ \hline
        \multicolumn{1}{|l||}{$10^2$}            & 86$^\dagger$ (2)  & 20   (2)        &86$^{\mathsection}$ (2) & 71 (2) \\ \hline
        \multicolumn{1}{|l||}{$10^4$}          & 86$^\dagger$ (3) & 31 (3) &        86$^{\mathsection}$ (2) & 71 (2) \\ \hline
        \multicolumn{1}{|l||}{$10^6$}          & 86$^\dagger$ (7) & 128 (7) & 86$^{\mathsection}$ (4) &  266 (6) \\ \hline
    \end{tabular}
    \caption{Total number of  the higher-dimensional problem solves required by Method~2 for the case
    study in Subsection \ref{subsec:example2}.  The number  of the Newton  iterations are in brackets. Within each case
        the construction of the multiscale flux basis is done only once, we mark by
        $\dagger$ (respectively ${\mathsection}$) common computations.}%
    \label{tab:results_case_newton}
\end{table}

\textbf{ For Method 2}, involving Newton's method for the linearization step,  the number of  higher-dimensional problem solves  are reported in Table~\ref{tab:results_case_newton}. As expected, Method 2 is more efficient
than Method~1 in terms of the  number of  higher-dimensional problem solves required, regardless of using multiscale flux
basis in the domain decomposition algorithm. Again, the number  of solves for the classical approach(DD) depends
on the used parameters. This table demonstrates (as shown with Method~1) that as
the value of $\beta_{\gamma}$ is increased, there is a point after which
Method~2 with multiscale flux basis is more efficient than without multiscale flux basis. In that case,  the gain  in the number of solves   becomes  more significant when decreasing  the value of $\alpha_{\gamma}$. Note that, in practice,
the simulations for Method~2 with multiscale flux basis are performed with negligible computational costs as we reused the  flux  basis inherited from
Method~1. This point together with the fact that the total number of solves required by  the entire simulation of \textit{case (i)}  is now reduced by $57\%$ as well as
that of \textit{case (i)} is reduced by $80\%$ showcase the performance of   Method~2 with multiscale flux basis.

%

To sum up, equipping Method~1 and 2 with  multiscale flux basis leads to   powerful tools  to solve complex  fracture network with important
savings in terms of the number of higher-dimensional problem solves. Note that, as  known,  one limitation of Method~2 involving Newton  method is  that a good initial value is usually required to obtain a solution.  A good combination of both methods can  also be used, in which  one can perform first some fixed-point iterations and then switch to Newton method. Concerning the computational costs, let us point out that the fixed-point algorithm of Method~1
requires at each iteration the assembly of the matrix corresponding to the linearization of the Darcy-Forchheimer equations and the solution of a linear system. The Newton method in Method~2 is slightly more expensive since
one has to assemble two matrices at each iteration and to update the right-hand side.


\subsection{Heterogeneous Forchheimer model}\label{subsec:example3}

In this example we assign to the two biggest fractures (one horizontal and
one vertical) high permeability while to the others low permeability. For the
highly permeable fractures we adopt the physical parameters of  \textit{case (i)}, while for those with lower permeabilities, the physical parameters corresponding to \textit{case (ii)} together with zero Forchheimer coefficient. In this case, we want to test the applicability of Method~1 with and without multiscale flux basis on highly  heterogeneous setting for both the permeability and
the flow models. We consider \textit{level 2} for the computation and, subsequently,
the local grids of the low permeable fractures are coarsened by having half of the
original number of elements resulting in 60 mortar unknowns instead of
84.

\begin{figure}[hbt]
    \centering
    \includegraphics[width=0.475\textwidth]{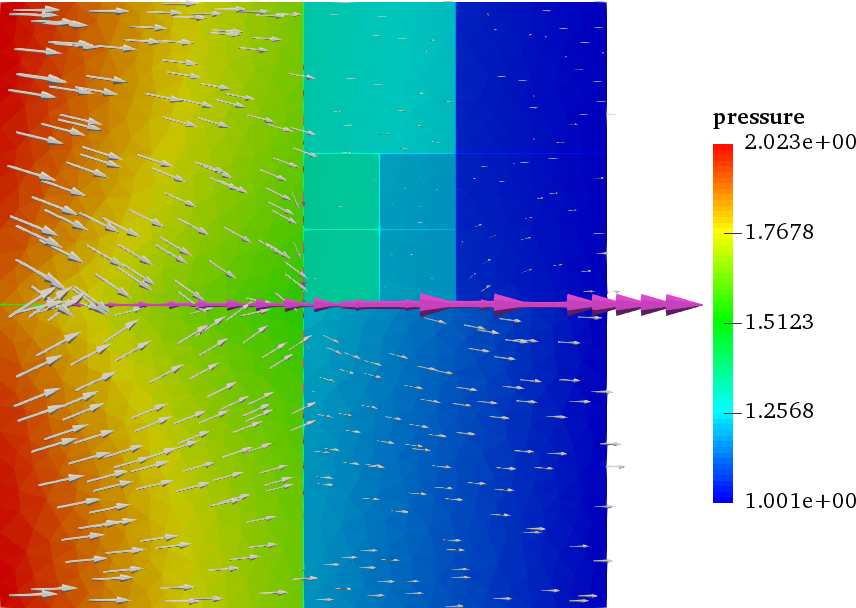}%
    \caption{Pressure and velocity solutions for example
    presented in Subsection \ref{subsec:example3} for $\beta = 10^2$.
    The velocity is represented by arrows (purple for the fractures)
    proportional to its magnitude.
    }%
    \label{fig:results_case3}
\end{figure}

As usually, we compare the method with and without multiscale flux basis in terms of the number of higher-dimensional problem solves. The results are represented in Table~\ref{tab:results_case3}.
\begin{table}[htbp]
    \centering
    \begin{tabular}{|l||l|l|}
        \hline
        $\beta_\gamma$  & MS      & DD    \\ \hline
        $1$             & 62$^\dagger$ (2)  & 63 (2) \\ \hline
        $10^2$          & 62$^\dagger$ (3)  & 84 (3) \\ \hline
        $10^4$          & 62$^\dagger$ (8)  & 189 (8)  \\ \hline
        $10^6$          & 62$^\dagger$ (64) & 2372 (64) \\ \hline
    \end{tabular}
    \caption{Total number of  the higher-dimensional problem solves required by Method~1 for the case
    study in Subsection \ref{subsec:example3}.  The number   of the fixed-point  iterations are in brackets. Within each case
    the construction of the multiscale basis is done only once, we mark by
    $\dagger$ common computations.}%
    \label{tab:results_case3}
\end{table}
In the present setting, we can see that the classical approach is outperformed with the approach equipped with multiscale flux basis,   particularly, the total computational costs is drastically reduced  when  the non-linear effects becomes more important. The entire simulation of Table~\ref{tab:results_case3} required 2708
higher-dimensional problem solves for the classical approach
while the same approach equipped with multiscale flux basis required 62 solves.  The overall gain is then of $94\%$ which
can also be appreciated for level 3. Similar
conclusions as above can be drawn for Method~2, namely in terms of reduction of the solves (not shown). An example of solution is given in Figure \ref{fig:results_case3}.

\subsection{Generalized Forchheimer model}\label{subsec:example4}
As stated previously, another  advantage distinguishes our approach
is that it  can integrate easily more complex problems. Here,  we apply our procedure to a more general model
describing the  pressure-flow relation in the fractures. Precisely, for
larger fracture flow velocities, the drag forces (in the Forchheimer model proportional
to the velocity norm) require to consider an additional term
proportional to the fluid viscosity. Considering the Barus formula
\cite{Barus1893}, we
have an exponential relation between the fluid viscosity and the pressure.
We  consider problem
\eqref{weak_mixed_formulation} where the non-linear term is as follows
\begin{gather*}
    \vecF^{-1}(\vecu_{\gamma},  p_\gamma) = \vK^{-1}_{\gamma}e^{\zeta
    p_\gamma}+\beta_{\gamma}\II|\vecu_{\gamma}|,
\end{gather*}
where $\zeta$ being a model parameter. Thus,  the non-linear effects are
now dependent on both   the pressure and  the velocity.
For a more detailed discussion we refer to \cite{Srinivasan2016}.
For the present setting, the  fracture permeabilities are set as in \textit{case (i)} and \textit{(iii)} of Table~\ref{tab:case1}.

For the discretization of the mixed geometry, we consider
\textit{level 2}. We use Method~1 with and without multiscale basis functions. Also, it was not necessary to recompute the basis functions,  since we can reuse
the stored multiscale flux basis from the previous test case and
solve then only on the fracture network the above more complex Darcy-Forchheimer model.
The total number  the higher-dimensional
problem solves  for $\beta_\gamma = 20$ and $\zeta \in (0.5, 5, 7.5)$ is
reported in Table~\ref{tab:results_case4}.
\begin{table}[htbp]
    \centering
    \begin{tabular}{l|l|l||l|l|}
        \cline{2-5}
        & \multicolumn{2}{c||}{\textit{case (i)}} & \multicolumn{2}{c|}{
        \textit{case (iii)}} \\
        \hline
        \multicolumn{1}{|l||}{$\zeta$} & MS      & DD       & MS & DD \\ \hline
        \multicolumn{1}{|l||}{$0.5$}   & 86$^\dagger$ (5)  &   71 (5) &
        86$^{\mathsection}$ (4) & 176 (4) \\ \hline
        \multicolumn{1}{|l||}{$5$}     & 86$^\dagger$ (4)  &  648 (4) &
        86$^{\mathsection}$ (6) & $\infty$ \\ \hline
        \multicolumn{1}{|l||}{$7.5$}    & 86$^\dagger$ (3)  & 5317 (3) &
        86$^{\mathsection}$ (4) &  $\infty$ \\ \hline
    \end{tabular}
    \caption{Total number of  the higher-dimensional problem solves required by Method~1 for the case     study in Subsection \ref{subsec:example4}.  The number   of the fixed-point  iterations are in brackets. Within each case
    the construction of the multiscale basis is done only once, we mark by
    $\dagger$ (respectively ${\mathsection}$) common computations.}%
    \label{tab:results_case4}
\end{table}
As expected, for such a strong non-linearity,  the results shows that
a considerable gain in terms of higher-dimensional problem solves
can be achieved. Particularly, for large values of
$\zeta$ the classical approach becomes uncompetitive to the new approach.
\begin{figure}[hbt]
    \centering
    \includegraphics[width=0.475\textwidth]{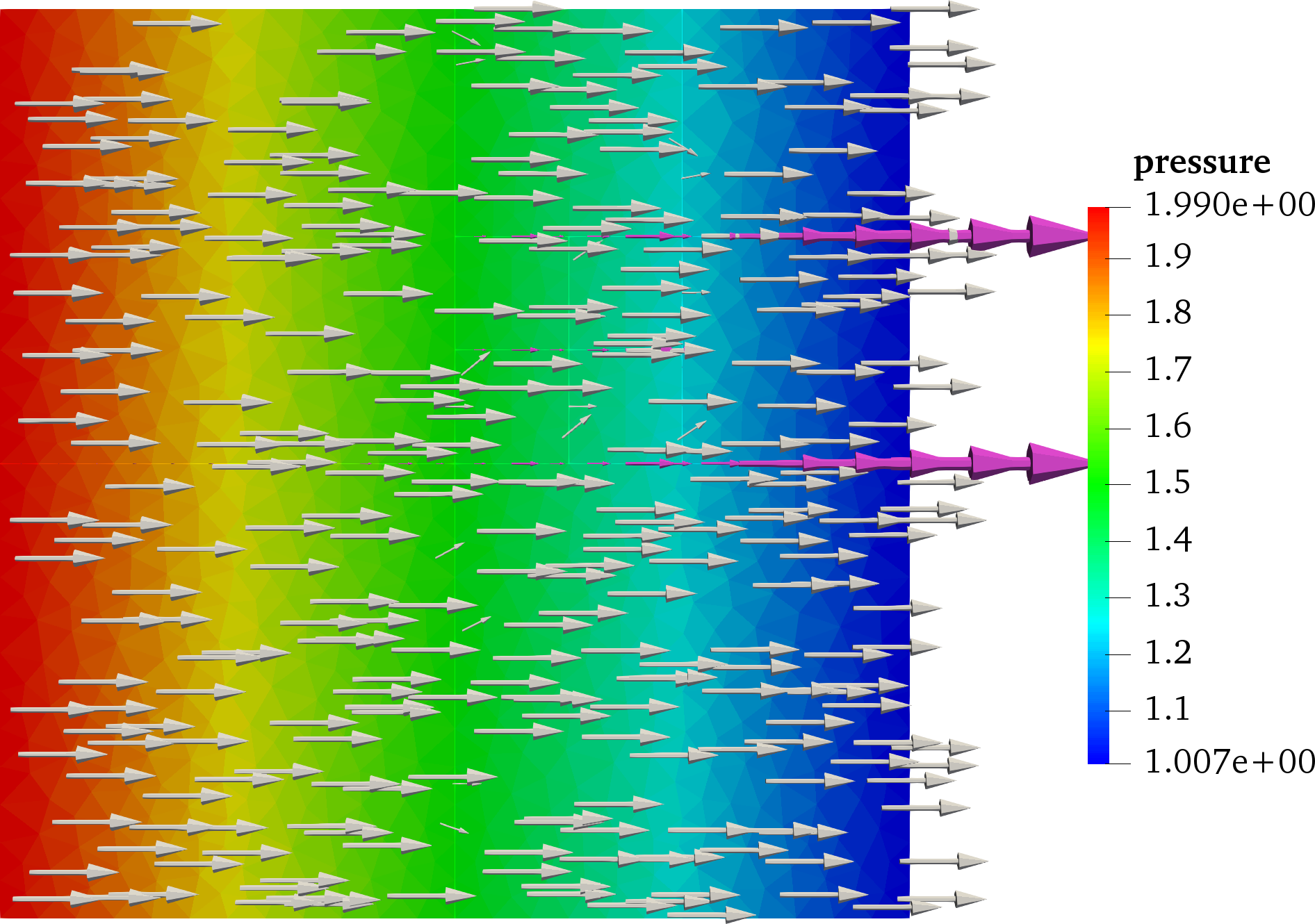}%
    \caption{Pressure and velocity solutions for example
    presented in Subsection \ref{subsec:example4} for $\zeta = 5$.
    The velocity is represented by arrows (purple for the fractures)
    proportional to its magnitude.
    }%
    \label{fig:results_case4}
\end{figure}
In Figure \ref{fig:results_case4} we report the solution for $\zeta = 5$.

\section{Conclusions}
\label{sec:conclusion}

In this work, we have presented a strategy to speed up the computation of a
Darcy-Forchheimer model for flow and pressure in fractured porous media by means
of multiscale flux basis, that represent the inter-dimensional flux exchange.
The scheme transforms a computationally expensive discrete fracture model to a
more affordable discrete fracture network, where in the latter only a
co-dimensional problem is solved. The multiscale flux basis are computed in an
\textit{offline} stage of the simulation and, despite the particular choice done
in this paper, are completely agnostic to the model in the fracture network. The
numerical results show the speed-up gain compared to a more classical
linearization--domain-decomposition approaches, where solves in both  the matrix
and the fracture network  are required along the entire outer--inner iterative
method. Crucially, an important number of the outer--inner interface
iterations may be spared.

With the proposed approach we are able to predict the computational effort
needed to solve the problem since it is directly related to the number of mortar
grids in the fracture network. Furthermore, the multiscale flux basis can be
reused when the fracture network geometry, rock matrix properties, and normal
permeability are fixed. Theoretical findings and numerical results show the
validity of the proposed approach and of its aforementioned properties.

Even if not explicitly considered in this work, it is possible to further
increase the efficiency of the proposed scheme by the following two steps.
First, compute a multiscale flux basis only in the related connected part of the
rock matrix.  Second, use an adaptive stopping criteria for the inner--outer iterative method based on a posteriori error estimates.  These enhancements are a part of future work along with
the extension in three-dimensions.

\section*{Acknowledgments}

We acknowledge the financial support from the Research Council of Norway for the TheMSES project (project no. 250223) and the ANIGMA project (project no. 244129/E20) through the ENERGIX program.
The authors also warmly thank Eirik Keilegavlen for valueable comments and discussions on this topic.

\bibliographystyle{siamplain}
\bibliography{biblio,biblioshear}
\end{document}